\documentclass[11pt]{article}
\usepackage{amssymb,latexsym,amsmath,graphicx}
\usepackage{mathtools}
\usepackage{latexsym}
\usepackage{mathrsfs}
\usepackage{amsfonts}
\usepackage{booktabs}

\DeclarePairedDelimiter\floor{\lfloor}{\rfloor}
\newcommand{\argmin}{\operatornamewithlimits{argmin}}

\newtheorem{thm}{Theorem}
\newtheorem{lem}{Lemma}
\newtheorem{rem}{Remark}
\newtheorem{cor}{Corollary}

\begin{document}

\begin{center}
	\Large{{\bf On regression and classification with possibly missing response variables in the data}}
\end{center}

\begin{center}
	{\rm Majid Mojirsheibani$^a$}
	\footnote{
		 Corresponding author.
 Email:  majid.mojirsheibani@csun.edu\\
		This work was supported by the National Science Foundation (NSF) under Grant DMS-1916161 of Majid Mojirsheibani},
{\rm William Pouliot$^b$, and 	Andre Shakhbandaryan$^c$}

\vspace{2mm}
	{\it $^{a,c}$ Department of Mathematics, California State University, Northridge, CA, USA
		
		$^b$ Department of Economics, University of Birmingham, Birmingham, UK}
\end{center}

\begin{abstract}
This paper considers the problem of kernel regression and classification with possibly unobservable response variables in the data, where the mechanism that causes the absence of information is unknown and can depend on both predictors and the response variables.  Our  proposed approach involves two steps: In the first step, we construct a family of models (possibly infinite dimensional) indexed by the unknown parameter of the missing probability mechanism.    In the second step, a search is carried out to find the empirically optimal member of an appropriate cover (or subclass) of  the underlying family in the sense of minimizing the mean squared prediction error.  The main focus of the paper is to look into the theoretical properties of these estimators. The issue of identifiability is also addressed. Our methods use a data-splitting  approach which is quite easy to implement.  We also derive exponential bounds on the performance of the resulting estimators in terms of their deviations from the true regression curve in general $L_p$ norms,  where we also allow the size of the cover or subclass to diverge as the sample size $n$ increases.  These bounds immediately yield various strong convergence results for the proposed estimators.  As an  application of our findings, we consider the problem of statistical classification based on the proposed regression estimators and also look into their rates of convergence under different settings.  
Although this work is mainly stated for kernel-type estimators, they can also be extended to other popular local-averaging methods such as nearest-neighbor estimators, and histogram estimators. 

\end{abstract}
{\bf MSC2020 subject classifications:} Primary 62G05; secondary 62G08

\vspace{3mm}\noindent
{\bf Keywords and phrases:} Regression, partially observed data, kernel,  convergence, classification, margin condition.

\allowdisplaybreaks

\section{Introduction}\label{sec.intro}
During the past decade, there has been a steady growing interest in developing appropriate procedures to perform estimation and inference in the presence of  incomplete data under the complex regime where the data is not missing at random (NMAR).  The NMAR setup is generally acknowledged to be a difficult problem in incomplete data literature due to identifiability issues; this is significantly different from the simpler  missing at random model where the absence of $Y$ depends on  ${\bf X}$ only (and not on both ${\bf X}$ and $Y$).   

\vspace{3mm}\noindent
The focus of this paper is on the important problem of the theoretical performance of kernel regression and classification under the realistic assumption that many response values in the data may be unavailable or missing.  Unobservable or incomplete  data occur frequently in medical data, survey data, public opinion polls,  as well as the data collected in many areas of scientific activities. Here, we assume that for various reasons some of the response values in the data may be unavailable or missing.  More specifically, let $({\bf X},Y )\in\mathbb{R}^d\times\mathbb{R}$ be a random  vector and consider the problem of estimating the regression function $m({\bf x}) = E(Y |{\bf X} = {\bf x})$, based on  $n$ independent and identically distributed (iid) observations $({\bf X}_i, Y_i)$,  $i = 1,\dots,n$, drawn from the distribution of $({\bf X}, Y)$. When the data is fully observable, the  classical Nadaraya-Watson kernel estimator of $m({\bf x})$ (Nadaraya (1964), Watson (1964))  is given by
\begin{equation}\label{D1}
\widehat{m}_n({\bf x}) =  \frac{\sum_{i=1}^n Y_i \,\mathcal{K}(({\bf x}-{\bf X}_i)/h_n)}{\sum_{i=1}^n \mathcal{K}({\bf x} - {\bf X}_i)/h_n)}\,,
\end{equation}
where the function $\mathcal{K}:\mathbb{R}^d\to\mathbb{R}_+$ is the kernel used with the bandwidth  $h_n>0$. A global measure of the accuracy of $\widehat{m}_n(\cdot)$, as an estimator of $m(\cdot)$, is given by its $L_{p}$-type statistic
$$
\begin{array}{cc}
I_n(p)=\displaystyle\int\left|\widehat{m}_n\left({\bf x}\right)-m\left({\bf x}\right)\right|^p\mu(d{\bf x}), & 1\leq p<\infty\,,
\end{array}$$
where $\mu$ is the probability measure of ${\bf X}$. The quantity $I_n(1)$ also plays an important role in statistical classification; see for example Devroye et al (1996, Sec. 6.2). For the strong convergence of $I_n(1)$ to zero see, for example,  Devroye and Krzy\`{z}ak (1989). In fact, in the cited paper, Devroye and Krzy\'{z}ak obtain a number of equivalent results under the assumption that $|Y|\leq L<\infty$, one of which states that if the kernel $\mathcal{K}$ is {\it regular} (the definition will be given later) then for every $\epsilon>0$ and $n$ large enough, one has
$P\left\{I_n(1)>\epsilon\right\}\leq \exp\{-c\,n\},$
where $c$ is a positive constant depending on $\epsilon$ but not  on $n$. 

\vspace{3.5mm}\noindent
Now, suppose that the response variable $Y$ is allowed to be missing according to the NMAR mechanism. Then, it is not hard to see that the estimator $\widehat{m}_n({\bf x})$ in (\ref{D1}) is no longer available. Of course, one might  decide (incorrectly) to use the kernel estimator based on the complete cases only, i.e., the estimator  
$m^{\mbox{\tiny cc}}_n({\bf x}) \,:= \,\sum_{i=1}^n \Delta_i Y_i \, 
\mathcal{K}(({\bf x}-{\bf X}_i)/h_n)\big/\sum_{i=1}^n \Delta_i \mathcal{K}(({\bf x}- {\bf X}_i)/h_n).$
Unfortunately, $m^{\mbox{\tiny cc}}_n({\bf x})$ turns out to be the estimator of the quantity $E(\Delta Y|X=x)\big/E(\Delta|{\bf X}={\bf x})$ which, in general, is not equal to the regression function $m({\bf x})$\,=\,$E(Y|{\bf X}={\bf x})$ under a NMAR response mechanism.  

\vspace{3.5mm}\noindent
For the important case of predictive models (such as regression), Kim and Yu (2011) considered a highly versatile logistic type missing probability
mechanism that works as follows. Define the indicator random variable  $\Delta$\,=\,0 if $Y$ is missing (and $\Delta$\,=\,1
otherwise), and let  $\pi({\bf x}, y) := E\left[\Delta\,\big| {\bf X}={\bf x}, Y=y\right]$ be the {\it selection probability}, also called the {\it nonresponse propensity}.  
Then, Kim and Yu (2011) considered the flexible model
\begin{equation}\label{NonIgnore}
\pi_{\gamma}({\bf x}, y) := E\left[\Delta\,\big| {\bf X}={\bf x}, Y=y\right]= \frac{1}{1+\exp\big\{g({\bf x})+\gamma y\big\}} \,,
\end{equation}
where $\gamma$ is an unknown parameter and $g$ is a completely unknown function of the predictor ${\bf x}$. In what follows, the true value of the unknown parameter $\gamma$ will be denoted by $\gamma^*$ everywhere. The missing probability mechanism (\ref{NonIgnore}) has been used and studied extensively in the literature; see, for example, 
Zhao and Shao (2015), Shao and Wang (2016), Morikawa et al (2017), Uehara and Kim (2018), 	Morikawa and Kim (2018), Morikawa and Kano (2018), Fang et al (2018), 	O'Brien et al (2018), 	Maity et al (2019), Sadinle and Reiter (2019), 	Zhao et al (2019), 	Yuan et al (2020), 	Chen et al (2020), 	Mojirsheibani (2021), and Liu and Yau (2021).  Of course, one may decide to consider more general nonparametric models instead of (\ref{NonIgnore}), but the estimation of such general models will become a difficult (if not impossible) issue. In fact, in view of the recent widespread use of model (\ref{NonIgnore}) in the literature, there appears to be the tacit consensus that (\ref{NonIgnore}) is versatile enough to be used in predictive models such as regression and classification, and this will also be the direction of the current paper. We observe that if  $\gamma$\,=\,0, then (\ref{NonIgnore}) reduces to the simpler case of missing at random assumption.

\vspace{3mm}\noindent
One of our aims in this paper is to explore the construction of several counterparts of the kernel estimator in (\ref{D1}) for the case where the response variable $Y$ can be missing. Another aim is to apply our results to the problem of classification where we construct asymptotically optimal nonparametric classification rules in the presence of NMAR data. Our contributions in these directions are three-fold. 
(i) We develop two types of easy-to-implement estimators of the regression curve $m({\bf x})$ in the presence of NMAR data.    Additionally, we consider a more general version of model (\ref{NonIgnore}) where the quantity $\exp\{\gamma y\}$ will be replaced by a more general positive function $\varphi(y)$.  We also propose  estimators of the component $\varphi(y)$ of our more general version of (\ref{NonIgnore}). 
The new estimators, which are based on the approximation theory of totally bounded classes of functions, are constructed using an easy-to-implement  data-splitting approach. 
(ii) We will  carefully explore and  study  the global properties of the proposed regression estimators in general $L_p$ norms; these results parallel those of Devroye and Krzy\`{z}ak (1989) for  the simpler case of no missing data. More specifically, we provide exponential performance bounds on the  $L_p$ norms of the proposed regression estimators that are valid under rather standard assumption.  Such bounds in conjunction with the Borel-Cantelli lemma immediately yield various strong convergence and optimality results. Exploiting these bounds further, we also look into the rates of convergence of the proposed estimators (in $L_p$). (iii)  A study of the applications of our proposed  estimators to the problem of nonparametric classification in the presence of partially observed data  is also considered.

\vspace{3mm}\noindent
As an important application of our results to the field of machine learning and statistical classification, we note that in the so-called semi-supervised learning one usually has to deal with large amounts of missing responses (or missing labels) in the data.  In such setups, researchers in machine learning have made efforts to develop procedures for utilizing the unlabeled cases (i.e., the data points with missing  $Y_i$'s) in order to construct more effective classification rules; see, for example, Wang and Shen (2007). But most such results assume that the response variable is missing completely at random; see, for example, Azizyan et al (2013). Our results in Section 3 makes it possible to develop classification rules in the presence of NMAR response variables for the semi-supervised setup, where we also study the rates of convergence of such classifiers.

\vspace{3mm}\noindent
The rest of the paper is organized as follows.  Section \ref{main} presents the main results, where in Subsection \ref{sub1} the estimation of the true $\gamma^*$ can be be based on any available method. It is shown in this case that the $L_p$-consistency of the resulting regression estimator requires the consistency of the estimator of $\gamma^*$ as well (which may be an issue). Subsection \ref{sub1} also proposes a generalization of the model (\ref{NonIgnore}), as given by (\ref{NonIg2}),  where new estimation methods based on the theory of totally bounded classes of functions are used to estimate the unknown function $\varphi^*$.  It is shown here that the resulting regression estimator can be $L_p$-consistent without requiring the consistency of the estimator of $\varphi^*$. Subsection \ref{HTTE} uses a Horvitz-Thompson type inverse weighting approach to estimate the underlying regression function.  Section \ref{ACPL} focuses on the applications of our estimators to the  problem of nonparametric classification with partially observed data. Here, we also look into the rates of convergence of the proposed classifiers under different conditions. All proofs are deferred to Section (\ref{PRF}).

\section{Main results}\label{main}
\subsection{The first estimator and a more general missing mechanism} \label{sub1}
Consider the missing probability mechanism (\ref{NonIgnore})  and let $\mathbb{D}_n=\{({\bf X}_1, Y_1, \Delta_1), \dots, ({\bf X}_n, Y_n, \Delta_n)\}$ be independent and identically distributed (iid) observations, i.e., the  data. Then, clearly  the estimator $\widehat{m}_n$ in (\ref{D1}) is no longer available due to the presence of missing $Y_i$'s. Furthermore, as discussed in the introduction, the complete-case estimator that only uses the fully observable data is not necessarily the correct estimator under model (\ref{NonIgnore}) anymore.  In the following two sections, we propose some alternative estimators instead.
To justify our first  estimator, we start by constructing an initial naive plug-in type estimator which works as follows.  Define the quantity
\begin{eqnarray}\label{ETA12}
\eta_k({\bf x}, t)=E\left[\Delta Y^{2-k}\exp\{t\, Y\}\Big| {\bf X}={\bf x}\right],~~~k=1, 2,~~~t\in\mathbb{R},
\end{eqnarray}
and observe that when $P(\Delta=1)\neq 1$ (i.e., when $Y$ is allowed to be missing), one can express the regression curve $m({\bf x})=E(Y|{\bf X}={\bf x})$ as 
\begin{equation} \label{repr}
m({\bf x})\,\equiv\,m_{\gamma^*}({\bf x}) \,=\, \eta_1({\bf x}, 0) +\frac{\eta_1({\bf x},\gamma^*)}{\eta_2({\bf x},\gamma^*)}\,\big(1-\eta_2({\bf x},0)\big),
\end{equation}
where $\gamma^*$ is the true value of $\gamma$ which is virtually always unknown; 
here, (\ref{repr}) follows from the more general representation of $m({\bf x})$ given by Lemma \ref{LEM-00} in this paper. Of course, $\gamma^*$ itself is unknown and has to be estimated.  Now, let 
$\widehat{\gamma}$ be any estimator of $\gamma^*$  (see Remark \ref{REM-AA} for identifiability issues)  and consider the following  simple kernel-type estimator of  (\ref{repr})
\begin{eqnarray}\label{mhat}
\widehat{m}_{n,\widehat{\gamma}}({\bf x}) &=& 
\widehat{\eta}_1({\bf x}, 0) +\frac{\widehat{\eta}_1({\bf x},\widehat{\gamma})}{\widehat{\eta}_2({\bf x},\widehat{\gamma})}\, 
\big(1-\widehat{\eta}_2({\bf x},0)\big),
\end{eqnarray}
where
\begin{eqnarray} \label{ETA12.hat}
\widehat{\eta}_k({\bf x}, t)=
\frac{\sum_{i=1}^n\Delta_i Y^{2-k}_i \exp\left\{t \,Y_i\right\}\mathcal{K}(({\bf x} -{\bf X}_i)/h)}{\sum_{i=1}^n\mathcal{K}(({\bf x}-{\bf X}_i)/h)},~~~
k=1, 2,~~~t\in\mathbb{R},
\end{eqnarray}
and, as in (\ref{D1}),   $\mathcal{K}:\mathbb{R}^d\to \mathbb{R}_+$ is the kernel used with bandwidth $h$. 
In passing, we also point out that although we are considering a kernel type estimator in (\ref{mhat}), virtually all our results in this paper continue to hold for other popular local-averaging estimators such as nearest-neighbor estimators, cubic histograms, as well as general  partitioning estimators. However,  to avoid making this work unnecessarily long and tedious, the paper is confined to kernel estimators only.

\vspace{4mm}\noindent
How good of an estimator is $\widehat{m}_{n,\widehat{\gamma}}({\bf x}) $ in (\ref{mhat})? To address and answer  this question, we start by assuming that the kernel $\mathcal{K}$ is {\it regular}: 

\vspace{2mm}\noindent
{\bf Definition} {\it A nonnegative kernel $\mathcal{K}$ is said to be regular  if there are real constants $b>0$ and $r>0$  such that $\mathcal{K}({\bf u})\geq b\, I\{{\bf u}\in S_{0,r}\}$ and $\int\sup_{{\bf y}\in {\bf u}+S_{0,r}} \mathcal{K}({\bf y})\,d{\bf u} < \infty$, where $S_{0,r}$  is the ball of radius $r$  centered at the origin.} 

\vspace{3mm}\noindent
For more on this, see  Devroye and Krzy\`{z}ak (1989). We also require the following condition regarding  the selection probability $\pi({\bf x}, y) := P(\Delta=1|{\bf X}={\bf x}, Y=y)$, which is quite standard in missing data literature: 

\vspace{4mm}\noindent
{\bf Assumption (A).} The selection probability,  $\pi({\bf x}, y)$, satisfies $\inf_{\boldsymbol{z},y}\, \pi(\boldsymbol{z}, y) =:\,\pi_{\mbox{\tiny min}}\,>\,0$,\, for some  $\pi_{\mbox{\tiny min}}$.

\vspace{3mm}\noindent
Assumption (A) essentially states that the response $Y$ can always be observed with a non-zero probability for any values of ${\bf x}$ and $y$. The following basic result  gives upper bounds on the performance of the $L_p$ norms  of the estimator $\widehat{m}_{n,\widehat{\gamma}}({\bf x}) $ under standard assumptions.

\begin{thm}\label{THM-A} 
Let $\widehat{m}_{n,\widehat{\gamma}}({\bf x}) $ be the estimator of $m({\bf x})$ defined in (\ref{mhat}), where $\widehat{\gamma}$ may be any  estimator of $\gamma^*$ in (\ref{repr}), and suppose that Assumption\,(A) holds. Suppose that the kernel $\mathcal{K}$ in (\ref{ETA12.hat}) is regular and that its bandwidth satisfies $h\to 0$ and $n h^d\to \infty$, as $n\to\infty$. Then, for every $\epsilon>0$, every $1\leq p <\infty$,  any distribution of $({\bf X}, Y)\in\mathbb{R}^d\times[-L, L]$, $L<\infty$, and $n$ large enough,
\begin{eqnarray}\label{Bound-1}
P\left\{\int\Big| \widehat{m}_{n,\widehat{\gamma}}({\bf x})   - m({\bf x})\Big|^p \mu(d{\bf x}) >\epsilon\right\}  &\leq& c_1\,e^{-c_2 n} 
+\, c_3\,P\big\{\left|\widehat{\gamma}-\gamma^*\right|> C_0\big\},
\end{eqnarray}
where $\mu$ is the probability measure of ${\bf X}$ and $c_1,\,c_2,\, c_3$, and  $C_0$ are positive constants not depending on $n$; here, $c_2$ also depends on $\epsilon$. 
\end{thm}
In passing, we note that  the bound in Theorem \ref{THM-A} is in the spirit of the classical result of Devroye and Krzy\`{z}ak (1989) for kernel regression estimators with no missing data (modulo the probability term $P\{|\widehat{\gamma}-\gamma^*|> C_0\}$
on the right side of (\ref{Bound-1})). Therefore, (\ref{Bound-1}) may be viewed as a generalization of the results of Devroye and Krzy\`{z}ak (1989), except that  we are allowing the response variable $Y$ to be missing not at random.

\begin{rem}\label{REM-AA}
The bound in Theorem \ref{THM-A} shows that the consistency of  $\widehat{\gamma}$, as an estimator of $\gamma^*$, is needed in order for the proposed regression estimator to converge in the $L_p$ norm. Unfortunately, due to parameter identifiability issues, consistent estimation of $\gamma^*$ can be a serious challenge unless one either has access to additional  external data, as in Kim and Yu (2011), or one can correctly assume that the function $g({\bf x})$ in  (\ref{NonIgnore}) is independent/free of certain components of ${\bf x}$\,=\,$(x_1,\cdots, x_d)^{\mbox{\tiny $T$}}$; see, for example, Shao and Wang (2016) or Uehara and Kim (2018).  Here, we consider a different estimation procedure based on the approximation theory of totally bound class of functions.
\end{rem}

\vspace{3mm}\noindent
In what follows,  we start by considering a more general version of the missing probability model  (\ref{NonIgnore}), given by
\begin{equation}\label{NonIg2}
\pi_{\varphi}({\bf x}, y) := E\left[\Delta\,\big| {\bf X}={\bf x}, Y=y\right]= 
\frac{1}{1+\exp\big\{g({\bf x})\big\}\cdot \varphi(y)}\,,
\end{equation}
for an unknown function  $\varphi(y)>0$; the true $\varphi$ will be denoted by $\varphi^*$.  Clearly the function $\exp\{\gamma y\}$ in (\ref{NonIgnore}) is a special case of $\varphi(y)$.  Our approach to estimate the function  $\varphi^*$ here is based on the approximation theory of totally bounded function spaces.  More specifically, consider the situation where $\varphi^*$  belongs to a totally bounded class of functions in the following sense:
Let   $\mathcal{F}$ be a given class of function $\varphi: [-L, L] \rightarrow (0,B],$ for some $B<\infty$.   Fix $\varepsilon>0$ and suppose that the finite collection of functions $\mathcal{F}_{\varepsilon}=
\{\varphi_1, \dots, \varphi_{\mbox{\tiny $N(\varepsilon)$}}\}$,  $\varphi_i:[-L, L] \to (0,B],$ is an $\varepsilon$-cover of $\mathcal{F}$, i.e., for each $\varphi\in \mathcal{F}$, there is a $\bar{\varphi}\in \mathcal{F}_{\varepsilon}$ such that $\lVert \varphi-\bar{\varphi}\rVert_{\infty} < \varepsilon$; here, $\parallel \parallel_{\infty}$ is the usual supnorm. The cardinality of the smallest $\varepsilon$-cover of $\mathcal{F}$  is called the {\it covering number} of the family $\mathcal{F}$ and will be denoted by $\mathcal{N}(\varepsilon, \mathcal{F})$. If $\mathcal{N}(\varepsilon, \mathcal{F})<\infty$ holds for every $\varepsilon>0$, then the family $\mathcal{F}$ is said to be {\it totally bounded} (with respect to $\parallel \parallel_{\infty}$).  The monograph by van der Vaart and Wellner (1996; p. 83) provides more details on such concepts.

\vspace{3.5mm}\noindent
To present our methods and results, we employ a data splitting approach that works as follows.  Let $\mathbb{D}_n=\{(\boldsymbol{X}_1,Y_1,\Delta_1),\dots,$  $((\boldsymbol{X}_n,Y_n,\Delta_n)\}$ represent the sample of size $n$ (iid), where $\Delta_i=0$ if $Y_i$ is missing (and $\Delta_i=1$ otherwise). Now, randomly split the data  into a training sample $\mathbb{D}_m$ of size $m$ and a validation sequence $\mathbb{D}_{\ell}$ of size $\ell=n-m$, where $\mathbb{D}_m\cup\mathbb{D}_{\ell} =\mathbb{D}_n$ and $\mathbb{D}_m\cap\mathbb{D}_{\ell} =\varnothing$. Here, it is assumed that $\ell\to\infty$ and $m\to \infty$, as $n\to \infty$; the choices of $m$ and $\ell$ will be discussed later in our main results. Also, define the index sets
\[
\boldsymbol{{\cal I}}_m=\Big\{i\in\{1,\cdots,n\}\,\Big|\,({\bf X}_i,Y_i,\Delta_i)\in \mathbb{D}_m\Big\}~~\mbox{and}~~~
\boldsymbol{{\cal I}}_\ell=\Big\{i\in\{1,\cdots,n\}\,\Big|\,({\bf X}_i,Y_i,\Delta_i)\in \mathbb{D}_\ell\Big\}.
\]
Next, for each fixed $\varphi\in\mathcal{F}$, consider  the kernel-type estimator of $m({\bf x})$ constructed based on the training set $\mathbb{D}_m$ alone, given by
\begin{eqnarray}\label{mhat3}
\widehat{m}_{m}({\bf x};\varphi) &=& 
\widehat{\eta}_{m,1}({\bf x}) +\frac{\widehat{\psi}_{m,1}({\bf x}; \varphi)}{\widehat{\psi}_{m,2}({\bf x}; \varphi)}\, 
\big(1-\widehat{\eta}_{m,2}({\bf x})\big),
\end{eqnarray}
where $\widehat{\psi}_{m,k}({\bf x}; \varphi)$ and $\widehat{\eta}_{m,k}({\bf x})$, $k=1,2$, are the quantities
\begin{eqnarray} 
\widehat{\psi}_{m,k}({\bf x};\varphi)&=&
\frac{\sum_{i\in\boldsymbol{{\cal I}}_m}\Delta_i Y^{2-k}_i \varphi(Y_i)\,\mathcal{K}(({\bf x} -{\bf X}_i)/h)}{\sum_{i\in\boldsymbol{{\cal I}}_m}\mathcal{K}(({\bf x}-{\bf X}_i)/h)},~~~k=1, 2,~~~\varphi\in\mathcal{F},\label{PSI12.hat}\\
\widehat{\eta}_{m,k}({\bf x}) &=&
\frac{\sum_{i\in\boldsymbol{{\cal I}_m}}\Delta_i Y^{2-k}_i \mathcal{K}(({\bf x} -{\bf X}_i)/h)}{\sum_{i\in\boldsymbol{{\cal I}_m}}\mathcal{K}(({\bf x}-{\bf X}_i)/h)},~~~k=1, 2. \label{ETA12m.hat}
\end{eqnarray}
Of course, (\ref{mhat3})  is not quite the right estimator because the true $\varphi^*$ is unknown. To estimate the function $\varphi^*$, we first observe that in view of the results of Kim and Yu (2011), the term $\exp\{g({\bf x})\}$ that appears in (\ref{NonIg2}) can  also be expressed as
\begin{equation}\label{expgx}
\exp\{g({\bf x})\} = E\big[1-\Delta\big|{\bf X}={\bf x}\big]\big/E\big[\Delta\, \varphi(Y)\big|{\bf X}={\bf x}\big].
\end{equation}
However, estimating the right side of (\ref{expgx}) can be challenging due to identifiability issues, and a sufficient condition for model identification is (see, for example, Uehara and Kim (2018)) to assume that there is a part of ${\bf X}$, say ${\bf V}$, which is conditionally independent of $\Delta$, given $Y$ and ${\bf Z}$, where ${\bf X}=({\bf Z}, {\bf V})$; see Assumption (G) on the next page. 
Under this assumption, the selection probability in (\ref{NonIg2}) becomes
\begin{equation}\label{NonIg4}
\pi_{\varphi}({\bf z}, y) := E\left[\Delta\,\big| {\bf Z}={\bf z}, Y=y\right]= \frac{1}{1+\exp\{g({\bf z})\}\cdot \varphi(y)}\,.
\end{equation}
It is not hard to see that under (\ref{NonIg4}) one has 
\begin{equation}\label{expgz2}
\exp\{g({\bf z})\} = E\big[1-\Delta\big|{\bf Z}={\bf z}\big]\big/E\big[\Delta\, \varphi(Y)\big|{\bf Z}={\bf z}\big].
\end{equation}
Next, we propose the following two-step procedure to estimate the unknown function $\varphi^*$: 

\vspace{2mm}\noindent
{\it Step 1.} For each fixed (given) $\varphi\in \mathcal{F}_{\varepsilon}$, and in view of (\ref{expgz2}), the selection probability  $\pi_{\varphi}({\bf z}, y)$ in (\ref{NonIg4}) is estimated, based on $\mathbb{D}_m$ alone, by 
\begin{eqnarray}\label{FIhat1}
\widehat{\pi}_{\varphi}({\bf z}, y)= 
\left[1\,+\, \widehat{\exp\{g({\bf z})\}} \cdot\varphi(y)\right]^{-1},
\end{eqnarray}
where in view of (\ref{expgz2}), $\widehat{\exp\{g({\bf z})\}}$ is given by the following kernel-type estimator
\begin{equation}\label{gz.hat}
\widehat{\exp\{g({\bf z})\}} = \frac{\sum_{i\in\boldsymbol{{\cal I}_m}}(1-\Delta_i)\mathcal{H}(({\bf z} - {\bf Z}_i)/\lambda)}{\sum_{i\in\boldsymbol{{\cal I}}_m}\Delta_i  \varphi(Y_i)\,\mathcal{H}(({\bf z} -{\bf Z}_i)/\lambda)}\,;
\end{equation}
here $\mathcal{H}$ is the kernel used with bandwidth $\lambda$.

\vspace{3mm}\noindent
{\it Step 2.}  Let $\varepsilon_n>0$ be a  decreasing sequence $\varepsilon_n$$\,\downarrow\,$0, as $n\to\infty$,
and let $\mathcal{F}_{\varepsilon_n}=\{\varphi_1,\dots, \varphi_{N(\varepsilon_n)}\}$  $\subset \mathcal{F}$ be any $\varepsilon_n$-cover of $\mathcal{F}$. The proposed estimator of $\varphi^*$ is then defined by 
\begin{equation}\label{SK1}
\widehat{\varphi}_{n}~:=~\argmin_{\varphi\,\in\mathcal{F}_{\varepsilon_n}}\,
\ell^{-1} \sum_{i\in\,\boldsymbol{{\cal I}}_\ell} \frac{\Delta_i}{ \widehat{\pi}_{\varphi}({\bf Z}_i, Y_i)}\,\big| 
\widehat{m}_{m}({\bf X}_i;\varphi)-Y_i\big|^2,
\end{equation}
where $\widehat{m}_{m}({\bf x};\varphi)$ is as in (\ref{mhat3}). The subscript $n$ at $\widehat{\varphi}_{n}$ reflects the fact that the entire data of size $n$ has been used here. Finally, the corresponding  estimator of the unknown regression function $m({\bf x})$ is given by
\begin{equation}\label{mhat5}
\widehat{m}({\bf x};\widehat{\varphi}_{n})\,:=~\widehat{m}_{m}({\bf x};\varphi)\big|_{\varphi =\widehat{\varphi}_{n}},~~~\mbox{with $\widehat{m}_{m}({\bf x};\varphi)$ as in (\ref{mhat3}).}
\end{equation}
The estimator  in (\ref{SK1}) may be viewed as the empirical version of the minimizer of the mean squared error, i.e., the empirical version of
\begin{equation}
\varphi_{\varepsilon_n} \,:=~ \argmin_{\varphi\in\mathcal{F}_{\varepsilon_n}}\, E\big|m({\bf X}; \varphi)-Y\big|^2, \label{FIEP}
\end{equation}
where $m({\bf X};\varphi)$  is the regression function $m({\bf X};\varphi^*)$ evaluated at an arbitrary $\varphi\in \mathcal{F}_{\varepsilon_n}$ (see Lemma \ref{LEM-00}).
We also note that $\varphi_{\varepsilon_n}$ in (\ref{FIEP}) is an approximation to the true function $\varphi^*$ based on the cover $\mathcal{F}_{\varepsilon_n}$ of $\mathcal{F}$. 
In fact, we have
\begin{equation}
\varphi^{*} \,:=~ \argmin_{\varphi:\,[-L,L]\to\mathbb{R}_+}\, E\big|m({\bf X}; \varphi)-Y\big|^2. \label{PHIstar}
\end{equation}
\begin{rem}\label{REM-HHH}
To appreciate that (\ref{SK1}) is the correct empirical version of $\varphi_{\varepsilon_n}$,  we observe that upon conditioning on both $Y$ and ${\bf X}$, one obtains
$
E\big\{\frac{\Delta}{\pi_{\varphi}({\bf Z}, Y)}\,
\big|m({\bf X};\varphi)-Y\big|^2\big\}=E\big[E\big\{
\frac{\Delta}{\pi_{\varphi}({\bf Z}, Y)}\,
\big|m({\bf X};\varphi)-Y\big|^2\big|{\bf X},Y\big\}\big]
=E\big|m({\bf X}; \varphi)-Y\big|^2,
$
where the last equality follows from the definition of  $\pi({\bf z},y)$   in (\ref{NonIg2}).  Of course, there are alternative ways to estimate $\varphi$; one could, for example, consider minimizing expressions such as $\sum_{i\in\,\boldsymbol{{\cal I}}_\ell} \left(\frac{\Delta_i}{ \widehat{\pi}_{\varphi}({\bf Z}_i, Y_i)}-1\right)^p\cdot W({\bf Z}_i, \varphi)$, for some specified weight function $W({\bf z}, y)$ and $p$=\,1\,or\,2. However, with such choices, we have not been able to study and track down the $L_p$-norms of the resulting kernel regression estimators. This is because our technical Lemma \ref{LEM-1} does not hold for such alternative estimators.
\end{rem}

\vspace{2mm}\noindent
How good is   $\widehat{m}({\bf x};\widehat{\varphi}_{n})$ 
in (\ref{mhat5}) as an estimator of the true regression curve $m({\bf x})$? To answer this, 
we first state a number of assumptions.

\vspace{3mm}\noindent
{\bf Assumption (B).}  The kernel $\mathcal{K}$ satisfies $\int_{\mathbb{R}^d} \mathcal{K}({\bf x})\,d{\bf x}=1$ and
$\int_{\mathbb{R}^d} |x_i| \mathcal{K}({\bf x})\,d{\bf x}<\infty,$ for  $x_i\in(x_1,\cdots,x_d)^{\mbox{\tiny $T$}}= {\bf x}$. Also, the smoothing
parameter $h$ satisfies $h\to 0$ and $m h^d\to \infty$, as
$n\to \infty$.

\vspace{3mm}\noindent
{\bf Assumption (C).} The density function $f({\bf z})$ of  ${\bf Z}$ is compactly supported  and is bounded away from zero and infinity on its compact support. Additionally, the first-order partial derivatives of $f$ exist and are bounded on the interior of its support.

\vspace{3mm}\noindent
{\bf Assumption (D).}  $E[\Delta\,\varphi(Y)|{\bf X}={\bf x}] \geq \varrho_0$, for $\mu$--a.e.\,${\bf x}$ and each $\varphi\in \mathcal{F},$ for some  constant $\varrho_0 >0$. 

\vspace{3mm}\noindent
{\bf Assumption (E).} The partial derivatives $\frac{\partial}{\partial z_i} E[\Delta|{\bf Z}={\bf z}]$ and $\frac{\partial}{\partial z_i} E[\Delta\,\varphi(Y)|{\bf Z}={\bf z}]$ exist for $i=1,\dots, \mbox{dim}({\bf z})$, and are bounded on the compact support of $f$.

\vspace{3mm}\noindent
{\bf Assumption (F).}  $\mathcal{F}$ is a totally bounded class of functions $\varphi: [-L, L] \to (0, B]$, for some $B<\infty$ and $L<\infty$. 

\vspace{3mm}\noindent
{\bf Assumption (G).} {\it [Identifiability]} There is a part of ${\bf X}$, say ${\bf V}$, which is conditionally independent of $\Delta$, given $Y$ and ${\bf Z}$, where ${\bf X}=({\bf Z}, {\bf V})$.

\vspace{3.5mm}\noindent
Assumption (B) is not restrictive at all because the choice of the kernel $\mathcal{K}$ is at our discretion. The first part of Assumption (C) is usually imposed in the literature on nonparametric regression to avoid unstable estimates of $m({\bf x})$ in the tails of the density, $f$. The second part of this assumption is technical. Assumption (D) is quite mild and is justified by the fact that $E[\Delta\,\varphi(Y)|{\bf X}] =
 E[\varphi(Y)E(\Delta|{\bf X}, Y)|{\bf X}]\geq \pi_{\min} E[\varphi(Y)|{\bf X}]$, together with the fact that $\varphi(y)>0$ for all $y$.  Assumption (E) is technical and has already been used in the literature. Assumption (G) is a standard sufficient condition for model identifiability; see, for example, Uehara and Kim (2018). 

\vspace{3.5mm}\noindent
The following result  gives exponential upper bounds on the performance of the $L_2$ norms  of the estimator defined via (\ref{mhat5}) in conjunction with (\ref{SK1}). This result readily extends to more general $L_p$ norms ($p\geq 2$); see Remark \ref{REM-thm2} below.
\begin{thm}\label{THM-B} 
Let $\widehat{m}({\bf x};\widehat{\varphi}_{n})$ be as in (\ref{mhat5}) and suppose that Assumptions (A)--(G) hold. Also let  the missing probability mechanism $\pi$ be as in (\ref{NonIg4}). Then for every $\varepsilon_n>0$ satisfying $\varepsilon_n \downarrow 0$, as $n\to\infty$, every $t>0$,  any distribution of $({\bf X}, Y)\in\mathbb{R}^d\times[-L, L]$, $L<\infty$, and $n$ large enough,
	\begin{equation}
		P\left\{\int\Big|\widehat{m}({\bf x}; \widehat{\varphi}_n)- m({\bf x}) \Big|^2  \mu(d{\bf x}) >  t \right\}
		\,\leq~ c_4|\mathcal{F}_{\varepsilon_n}| e^{-c_5\ell t^{2}} 
		+ c_6 \ell\, |\mathcal{F}_{\varepsilon_n}|\left(e^{-c_7 m h^d}+ e^{-c_{8} m h^d\, t^{2}} \right),\label{THM2i}
	\end{equation}
	whenever $\varphi^*\in\mathcal{F}$, where  $|\mathcal{F}_{\varepsilon_n}|$ is the cardinality of the set 
	$\mathcal{F}_{\varepsilon_n}$ and $c_4$--\,$c_{8}$ are positive constants not depending on $m$, $\ell$, $n$, or $t$.
\end{thm}
\begin{rem} \label{REM-thm2}
Although the above theorem is stated in the $L_2$ sense, the theorem continues to hold for all $p\geq 2$. To appreciate this, observe that in the case of  $p>2$ one can always write
\[
\big|\widehat{m}({\bf x}; \widehat{\varphi}_n)- m({\bf x}) \big|^p \leq \Big(\big|\widehat{m}({\bf x}; \widehat{\varphi}_n)\big| + \big|m({\bf x}) \big|\Big)^{p-2}\big|\widehat{m}({\bf x}; \widehat{\varphi}_n) - m({\bf x}) \big|^2\leq (3L)^{p-2}\big|\widehat{m}({\bf x}; \widehat{\varphi}_n)- m({\bf x}) \big|^2.
\]
On the other hand, if $p\in[1, 2)$ then by H\"{o}lder's inequality we have
\[
P\left\{\int\Big|\widehat{m}({\bf x}; \widehat{\varphi}_n)- m({\bf x}) \Big|^p  \mu(d{\bf x}) > t\right\} ~\leq~ 
P\left\{\int\Big|\widehat{m}({\bf x}; \widehat{\varphi}_n)- m({\bf x}) \Big|^2  \mu(d{\bf x}) > t^{2/p}\right\}.
\]
\end{rem}
\vspace{3mm}\noindent
In passing, we note that the  bound in (\ref{THM2i}) may be viewed as a generalization of the classical results of Devroye and Krzyzak (1989) for kernel regression estimators with fully observable data.  Furthermore, as the following simple corollary shows, the above theorem can also be used to establish strong convergence  results.
\begin{cor}\label{COR-A}
Let $\widehat{m}({\bf X}; \widehat{\varphi}_n)$ be  the  estimator in (\ref{mhat5}). If, as $n\to\infty$, 
\begin{equation} \label{cond}
\varepsilon_n\downarrow 0,~~~\frac{\log \ell}{m h^d} \to 0, ~~~ \frac{\log|\mathcal{F}_{\varepsilon_n}|}{m h^d}\to 0,~\mbox{and}~~~\frac{\log|\mathcal{F}_{\varepsilon_n}|}{\ell}\to 0,
\end{equation}
then, under the conditions of Theorem \ref{THM-B},  we have
\[
E\Big[\big|\widehat{m}({\bf X}; \widehat{\varphi}_n)- m({\bf X}) \big|^p\Big|\mathbb{D}_n\Big]\longrightarrow^{a.s.} 0,~~~ \mbox{for all $p\in [2,\infty)$.}
\]
\end{cor}

\vspace{2mm}\noindent
Clearly, by Lebesgue dominated convergence theorem,  under the conditions of Corollary \ref{COR-A}, and without further ado, 
\[
E\big|\widehat{m}({\bf X}; \widehat{\varphi}_n)- m({\bf X}) \big|^p\longrightarrow 0,~~~ \mbox{for all $p\in [2,\infty)$.}
\]
Unfortunately, this result does not provide a rate of convergence.  The following theorem sheds more light on the convergence properties of the estimator in (\ref{mhat5}).

\begin{thm}\label{THM-BB} 
	Consider the estimator  $\widehat{m}({\bf X}; \widehat{\varphi}_n)$ in (\ref{mhat5}).  Then, under the  conditions of Theorem \ref{THM-B}, for $n$ large enough,
	\begin{eqnarray*}
	&& E\Big|\widehat{m}({\bf X}; \widehat{\varphi}_n)- m({\bf X}) \Big|^p \\ [3pt]
	&&~~~\leq~ \sqrt{\frac{a_1+\log \ell+ \log|\mathcal{F}_{\varepsilon_n}|}{a_2\cdot(\ell \wedge mh^d)}}
	\,+ \sqrt{\frac{1}{a_3\cdot(\ell \wedge mh^d)\big[ a_1+\log \ell+\log|\mathcal{F}_{\varepsilon_n}|\big]}}
	\,\,+ a_4\,\big|\mathcal{F}_{\varepsilon_n}\big|\,\ell\, e^{-c_7\,mh^d},~~
	\end{eqnarray*}
	for all $p\in [2,\infty)$, where $a_1$--\,$a_4$ are positive constants not depending on $m$, $\ell$, or $n$.
\end{thm}
The following result, which is an 
immediate corollary to Theorem \ref{THM-BB}, looks into the rate of convergence of the proposed regression estimator.
\begin{cor}\label{COR-A2}
	Let  $\widehat{m}({\bf X}; \widehat{\varphi}_n)$ be the estimator in (\ref{mhat5}) and suppose that (\ref{cond}) holds.  Then, under the  conditions of Theorem \ref{THM-B}, for all $p\geq 2$,
	\begin{equation*} 
	E\Big|\widehat{m}({\bf X}; \widehat{\varphi}_n)- m({\bf X}) \Big|^p \,=\, \mathcal{O}\left( \sqrt{\frac{\log(\ell \vee |\mathcal{F}_{\varepsilon_n}|)}{\ell \wedge mh^d}}\right).
	\end{equation*}
	In the special case where $m=\alpha\cdot n$ and $\ell=(1-\alpha)\cdot n$, where $\alpha\in(0,1)$, one finds (under the above conditions) that for all $p\geq 2$,
	\begin{equation*} 
	E\Big|\widehat{m}({\bf X}; \widehat{\varphi}_n)- m({\bf X}) \Big|^p \,=\, \mathcal{O}\left( \sqrt{\frac{\log(n \vee |\mathcal{F}_{\varepsilon_n}|)}{n h^d}}\right).
	\end{equation*}
\end{cor}

\vspace{3.5mm}\noindent
{\bf An Example.}\\  To compare and contrast the asymptotic performance of our estimation  approach with the existing methods, consider the class $\mathcal{F}$ of functions $\varphi$ of the form:
\begin{equation}\label{Exam-1}
\varphi(y)=\exp\{\gamma y\}\,,~~|\gamma|\leq M,~\, |y|\leq L,~~\mbox{for some  $M, L<\infty$},
\end{equation}
which is similar to the selection probability model used by Kim and Yu (2011). Clearly, if a consistent estimator $\widehat{\gamma}$ of $\gamma^*$ (the true value of $\gamma$) is available, then under the conditions of Theorem\,\ref{THM-A}, one immediately obtains the $L_p$ consistency of the corresponding regression function $\widehat{m}_{n, \widehat{\gamma}}({\bf x})$ in (\ref{mhat}). Unfortunately, such a consistent estimator is not readily available; for example, the estimator of Kim and Yu (2011) requires access to some  external data. Moreover, to use  the estimator of Shao and Wang (2016), it is necessary to be able to find a part ${\bf v}$ of the vector ${\bf x} = ({\bf u}, {\bf v})$ which is not involved in the function $g$  in  (\ref{NonIgnore}), i.e., one has to work with some $g({\bf u})$ instead of $g({\bf x})$ in  (\ref{NonIgnore}). Our estimator, however, evades such requirements in the sense that it bypasses the direct consistent estimation of the true value of $\gamma$.  In fact, for the function $\varphi$ as in (\ref{Exam-1}), it is straightforward to see that for every $\varepsilon>0$, the finite collection of functions
\begin{equation}\label{eps-cover}
	\mathcal{F}_{\varepsilon} = \left\{  \exp\{\gamma y \},~ |y|\leq L \bigg|\,
	\gamma \in \Big\{\left\{
	2\,i \varepsilon/\big(L\exp(ML)\big)\,\Big|\,|i|\leq \floor*{ML\exp\{ML\}/\varepsilon}\right\}\cup \{-M\}\cup\{M\}\Big\}
	\right\}
\end{equation}
is an $\varepsilon$-cover of $\mathcal{F}$ and the covering number of $\mathcal{F}$ is bounded by
$(2\,ML\exp\{ML\}\varepsilon^{-1}+3)$; see the Appendix for details. Since this bound grows like $\varepsilon^{-1}$ (as $\varepsilon\downarrow 0$), one obtains strong $L_p$ consistency results for the regression estimator (\ref{mhat5}) under the conditions of Theorem \ref{THM-B} for any sequence $\varepsilon_n\downarrow 0$ (as $n\to\infty$) for which $\log(1/\varepsilon_n)/(m h^d \vee \ell))\to 0.$ Similarly, the conclusions of Theorem \ref{THM-BB} and Corollary $\ref{COR-A2}$ continue to hold for such a sequence.

\subsection{A Horvitz-Thompson type estimator}\label{HTTE}
Our estimators in this section are based on a Horvitz-Thompson type inverse weighting approach (Horvitz and Thompson (1952)). This method works by scaling each observed response variable $Y$ with the inverse of the estimate of the {\it selection} probability, $\pi_{\varphi^*}({\bf Z},Y)$, as given by (\ref{NonIg4}), where $\varphi^*$ is the true function $\varphi$ in (\ref{NonIg4}) in the sense that
\begin{equation}
m({\bf X}; \pi_{\varphi^*}):=E\big[\Delta Y/\pi_{\varphi^*}({\bf Z},Y)\big|{\bf X}\big]=E[Y|{\bf X}]=m({\bf X}). \label{FS15}
\end{equation}
To motivate and describe this approach formally, 
consider the hypothetical (and unrealistic) situation where the true function $\pi_{\varphi^*}$ is completely known. Then a  kernel-type estimator of the regression curve $m({\bf x})$ based on inverse weighting is simply
\begin{equation}\label{BQ1}
\widetilde{m}_{n}({\bf x}; \pi_{\varphi^*})= \sum_{i=1}^n \, \frac{\Delta_i Y_i}{\pi_{\varphi^*}({\bf Z}_i, Y_i)}\, \mathcal{K}\big(({\bf x}-{\bf X}_i)/h\big)\Big/
 \sum_{i=1}^n  \mathcal{K}\big(({\bf x}-{\bf X}_i)/h\big).
 \end{equation}
Since $\pi_{\varphi^*}$ is unknown, we proceed as follows.  For each $\varphi\in {\cal F}$, consider the estimate of the selection probability $\pi_{\varphi}$ of (\ref{NonIg4}), based on $\mathbb{D}_m$, given by 
\begin{eqnarray}\label{CEQ1}
\widetilde{\pi}_{\varphi}({\bf Z}_i, Y_i)= 
\left[1\,+\, \frac{1-\widetilde{\eta}_{m}({\bf Z}_i)}{\widetilde{\psi}_{m}({\bf Z}_i; \varphi)} \cdot\varphi(Y_i)\right]^{-1},
\end{eqnarray}
where 
\begin{eqnarray} \label{ETPS}
\left\{
\begin{array}{ll}
\widetilde{\psi}_{m}({\bf Z}_i;\varphi) \,=\,
\sum_{j\in\boldsymbol{{\cal I}}_m ,\, j\neq i}\Delta_j\varphi(Y_j)\, \mathcal{H}(({\bf Z}_i -{\bf Z}_j)/h)\Big/\sum_{j\in\boldsymbol{{\cal I}}_m,\, j\neq i}\mathcal{H}(({\bf Z}_i-{\bf Z}_j)/h) ~& \\[8pt]
\widetilde{\eta}_{m}({\bf Z}_i) \,=\,
\sum_{j\in\boldsymbol{{\cal I}_m},\, j\neq i}\Delta_j  \mathcal{H}(({\bf Z}_i -{\bf Z}_j)/h)\Big/\sum_{j\in\boldsymbol{{\cal I}_m},\,j\neq i}\mathcal{H}(({\bf Z}_i-{\bf Z}_j)/h). ~ &
\end{array}
\right. 
\end{eqnarray}
Since $\pi_{\varphi}>\pi_{\min}>0 $ (by Assumption (A)) and since $\widehat{\psi}_{m}({\bf Z}_i;\varphi)$ in (\ref{ETPS}) is the estimator of $E[\Delta_i\varphi(Y_i)|{\bf Z}_i]\geq \varrho_0>0$, almost surely (see Assumption (D), we also consider the following truncated-type  version of the estimator in (\ref{CEQ1})
\begin{eqnarray}\label{CEQ2}
\breve{\pi}_{\varphi}({\bf Z}_i, Y_i)= 
\left[1\,+\, \frac{1-\widetilde{\eta}_{m}({\bf Z}_i)}{\pi_0\vee \widetilde{\psi}_{m}({\bf Z}_i; \varphi)} \cdot\varphi(Y_i)\right]^{-1},
\end{eqnarray}
where $\pi_0>0$ is a fixed constant whose choice will be discussed later  under Assumption (A$'$). Here, we note that $\breve{\pi}_{\varphi}$ in (\ref{CEQ2}) may be viewed as a one-sided winsorized estimator of $\pi_{\varphi}$ (compare this with $\widetilde{\pi}_{\varphi}$ in (\ref{CEQ1})). In applications with either simulated or real data, $\pi_0$ is chosen to be a small positive number such as $10^{-\nu},$ $\nu\geq 3$. Next,  let $\varepsilon_n>0$ be a  decreasing sequence $\varepsilon_n$$\,\downarrow\,$0, as $n\to\infty$
and let $\mathcal{F}_{\varepsilon_n}=\{\varphi_1,\dots, \varphi_{N(\varepsilon_n)}\}$  $\subset \mathcal{F}$ be any $\varepsilon_n$-cover of $\mathcal{F}$. Then, depending on whether (\ref{CEQ1}) or (\ref{CEQ2}) is used, an estimator of the unknown function $\varphi^*$ based on the $\varepsilon_n$-cover  $\mathcal{F}_{\varepsilon_n}$ is given by  
\begin{equation}\label{CEQ3}
\left\{
\begin{array}{ll}
\widetilde{\varphi}_{n}~:=~\argmin_{\varphi\,\in \mathcal{F}_{\varepsilon_n}}\, \ell^{-1} \sum_{i\in\, \boldsymbol{{\cal I}}_\ell} \frac{\Delta_i}{ \widetilde{\pi}_{\varphi}({\bf Z}_i, Y_i)}\,\big| 
\widehat{m}^{\mbox{\tiny HT}}_m({\bf X}_i;\widetilde{\pi}_{\varphi})-Y_i\big|^2,
~ & ~ \hbox{if (\ref{CEQ1}) is used,}  \\[8pt]
\breve{\varphi}_{n} ~:=~\argmin_{\varphi\,\in\mathcal{F}_{\varepsilon_n}}\,
\ell^{-1} \sum_{i\in\,\boldsymbol{{\cal I}}_\ell} \frac{\Delta_i}{ \breve{\pi}_{\varphi}({\bf Z}_i, Y_i)}\,\big| 
\widehat{m}^{\mbox{\tiny HT}}_m({\bf X}_i;\breve{\pi}_{\varphi})-Y_i\big|^2,
~ & ~ \hbox{if (\ref{CEQ2}) is used,}
\end{array}
\right. 
\end{equation}
where
\begin{eqnarray}\label{CEQ4}
\widehat{m}^{\mbox{\tiny HT}}_m({\bf x};\widetilde{\pi}_{\varphi}) &=& \sum_{i\in\boldsymbol{{\cal I}}_m}\, \frac{\Delta_i Y_i}{\widetilde{\pi}_{\varphi}({\bf Z}_i, Y_i)}\, \mathcal{K}\big(({\bf x}-{\bf X}_i)/h\big)\Big/  \sum_{i\in\boldsymbol{{\cal I}}_m}  \mathcal{K}\big(({\bf x}-{\bf X}_i)/h\big), 
\end{eqnarray}
and  $\widehat{m}^{\mbox{\tiny HT}}_m({\bf X}_i;\breve{\pi}_{\varphi})$ is obtained by replacing $\widetilde{\pi}_{\varphi}$ with $\breve{\pi}_{\varphi}$ in (\ref{CEQ4}). Finally, our proposed Horvitz-Thompson type estimator of the regression function $m({\bf x})$ is given by
\begin{equation} \label{CEQ5}
\left\{
\begin{array}{ll} 
\widehat{m}^{\mbox{\tiny HT}}({\bf x}; \widetilde{\pi}_{\widetilde{\varphi}_n})\,:=~\widehat{m}^{\mbox{\tiny HT}}_m({\bf x};\pi_{\varphi})\Big|_{\pi_{\varphi} =\widetilde{\pi}_{\widetilde{\varphi}_{n}}}
~ & ~ \hbox{if (\ref{CEQ1}) is used,} \\[8pt]
\widehat{m}^{\mbox{\tiny HT}}({\bf x}; \breve{\pi}_{\breve{\varphi}_n})\,:=~\widehat{m}^{\mbox{\tiny HT}}_m({\bf x};\pi_{\varphi})\Big|_{\pi_{\varphi} =\breve{\pi}_{\breve{\varphi}_{n}}}
~ & ~ \hbox{if (\ref{CEQ2}) is used,}
\end{array}
\right. 
\end{equation}
where $\widehat{m}^{\mbox{\tiny HT}}_m({\bf x};\pi_{\varphi})$ is as in (\ref{CEQ4}) but with $\widetilde{\pi}_{\varphi}$ replaced by $\pi_{\varphi}$. 

\vspace{3.5mm}\noindent
Next, we compare and study the asymptotic performance of the two estimators in (\ref{CEQ5}). It turns out, as in Theorem \ref{THM-B} and its corollary (i.e., Corollary \ref{COR-A}), that exponential upper bounds along with strong consistency results are available for both estimators. However, in the case of the winsorized-type estimator  $\widehat{m}^{\mbox{\tiny HT}}({\bf x}; \breve{\pi}_{\breve{\varphi}_n})$, one can also study the rates of convergence.   To state our results here, we first state the following counterpart of Assumption (A):

\vspace{4mm}\noindent
{\bf Assumption (A$'$).} The selection probability,  $\pi({\bf z}, y)$, satisfies $\inf_{{\bf z},y}\, \pi({\bf z}, y) =:\,\pi_{\mbox{\tiny min}}\,>\,0$ for some $\pi_{\mbox{\tiny min}}$ and the truncation constant $\pi_0$ in (\ref{CEQ2}) is any constant satisfying $0<\pi_0\leq \pi_{\mbox{\tiny min}}$.

\begin{thm}\label{THM-BBC} 
	Consider 
	the two regression function estimators defined via (\ref{CEQ5}) and let  the missing probability mechanism $\pi$ be as in (\ref{NonIg4}). 
	
	\vspace{2mm}\noindent
	\textbf{\textit{(i)}}~ Let $\widehat{m}^{\mbox{\tiny \rm HT}
		}({\bf x}; \widetilde{\pi}_{\widetilde{\varphi}_n})$ be the top estimator in (\ref{CEQ5}) and suppose that assumptions (A)--(G) hold. Then for every $\varepsilon_n>0$ satisfying $\varepsilon_n \downarrow 0$, as $n\to\infty$, every $t>0$,  any distribution of $({\bf X}, Y)\in\mathbb{R}^d\times[-L, L]$, $L<\infty$, and $n$ large enough,
	\begin{equation}\label{THM4i}
	P\left\{\int\Big| \widehat{m}^{\mbox{\tiny \rm HT}}({\bf x}; \widetilde{\pi}_{\widetilde{\varphi}_n}) - m({\bf x}) \Big|^2  \mu(d{\bf x}) \,>\,  t \right\}
	\,\leq~ \big|\mathcal{F}_{\varepsilon_n}\big| \,\Big( c_{9}\,e^{-c_{10} \ell t^2} + c_{11} \ell m \,e^{-c_{12} m h^d t^2}+ c_{13}\, \ell m \,e^{-c_{14}m h^d} \Big),
	\end{equation}
	whenever $\varphi^*\in\mathcal{F}$, where  $|\mathcal{F}_{\varepsilon_n}|$ is the cardinality of the set 
	$\mathcal{F}_{\varepsilon_n}$ and $c_9$--\,$c_{14}$ are positive constants not depending on $m$, $\ell$, $n$, or $t$.
	
	\vspace{2mm}\noindent
	\textbf{\textit{(ii)}}~ Let $\widehat{m}^{\mbox{\tiny\rm HT}}({\bf x}; \breve{\pi}_{\breve{\varphi}_n})$ be the second estimator in (\ref{CEQ5}) and suppose that assumptions (A$'$), (B)--(G) hold.  Then, under the conditions of part \textbf{\textit{(i)}} of the theorem, the bound in (\ref{THM4i}) continues to hold (with different constants $c_9$--$c_{14}$) for the probability $P\left\{\int |\widehat{m}^{\mbox{\tiny \rm HT}}({\bf x}; \breve{\pi}_{\breve{\varphi}_n})- m({\bf x}) |^2  \mu(d{\bf x}) >  t \right\}$. 
\end{thm}
\begin{rem} \label{REM-thm4}
	As in Remark \ref{REM-thm2}, it is straightforward to show that Part (ii) of the above theorem holds more generally for all $p\geq 2$. That is, the bound in (\ref{THM4i}) holds for $P\big\{\int |\widehat{m}^{\mbox{\tiny HT}}({\bf x}; \breve{\pi}_{\breve{\varphi}_n})- m({\bf x}) |^p  \mu(d{\bf x})$ $ >  t \big\}$, for all $p\geq 2$.
\end{rem}

\vspace{3.5mm}\noindent
The following result may be viewed as the counterpart of Corollary \ref{COR-A} for the two regression function estimators in (\ref{CEQ5}).
\begin{cor}\label{COR-AB}
	Consider the two estimators in (\ref{CEQ5}). If, as $n\to\infty$, 
	\begin{equation} \label{cond2}
	\varepsilon_n\downarrow 0,~~~\frac{\log (m \vee \ell)}{m h^d} \to 0, ~~~ \frac{\log|\mathcal{F}_{\varepsilon_n}|}{m h^d}\to 0,~\mbox{and}~~~\frac{\log|\mathcal{F}_{\varepsilon_n}|}{\ell}\to 0,
	\end{equation}
	then, under the conditions of Theorem \ref{THM-BBC},  the top estimator in (\ref{CEQ5}) satisfies the strong convergence property, $E\Big[\big|\widehat{m}^{\mbox{\tiny \rm HT}}({\bf X}; \widetilde{\pi}_{\widetilde{\varphi}_n})- m({\bf X}) \big|^2\Big|\mathbb{D}_n\Big]\rightarrow^{a.s.} 0.$ However, for the second estimator in (\ref{CEQ5}), 
	\begin{eqnarray*}
	&& E\Big[\big|\widehat{m}^{\mbox{\tiny \rm HT}}({\bf X}; \breve{\pi}_{\breve{\varphi}_n})- m({\bf X}) \big|^p\Big|\mathbb{D}_n\Big]\longrightarrow^{a.s.} 0\,,~~\mbox{for all $p\geq 2$}.
	\end{eqnarray*}
\end{cor}

\vspace{2mm}\noindent
In passing, we also note that  under the conditions of Corollary \ref{COR-AB}, by Lebesgue dominated convergence theorem, and without further ado, one has~
$
E\big|\widehat{m}^{\mbox{\tiny HT}}({\bf X}; \breve{\pi}_{\breve{\varphi}_n})- m({\bf X}) \big|^p\rightarrow 0,~ \mbox{for all $p\in [2,\infty)$.}
$
However, to study the rates of convergence here, we state the following theorem which is the counterpart of Theorem \ref{THM-BB} for the estimator $\widehat{m}^{\mbox{\tiny HT}}({\bf x}; \breve{\pi}_{\breve{\varphi}_n})$.

\begin{thm}\label{THM-BBCC} 
	Let  $\widehat{m}^{\mbox{\tiny HT}}({\bf x}; \breve{\pi}_{\breve{\varphi}_n})$ be the second estimator in (\ref{CEQ5}).  Then, under the  conditions of Theorem \ref{THM-BBC}, for all $p\in [2,\infty)$ and $n$ large enough,
	\begin{eqnarray*}
		&& E\Big|\widehat{m}^{\mbox{\tiny HT}}({\bf X}; \breve{\pi}_{\breve{\varphi}_n})- m({\bf X}) \Big|^p \\ [3pt]
		&&~~~\leq~ \sqrt{\frac{a_5+\log \ell+\log m + \log|\mathcal{F}_{\varepsilon_n}|}{a_6\cdot(\ell \wedge mh^d)}} \,+ \sqrt{\frac{1}{a_7 \cdot(\ell \wedge mh^d)\big[ a_5+\log \ell+\log m + \log|\mathcal{F}_{\varepsilon_n}|\big]}}~~~\\
		&&~~~~~~~~~~~+~a_8\,\big|\mathcal{F}_{\varepsilon_n}
		\big|\,\ell\,m\, e^{-c_7\,mh^d},
	\end{eqnarray*}
	where $a_5$--\,$a_8$ are positive constants not depending on $m$, $\ell$, or $n$.
\end{thm}
The following is an 
immediate corollary to Theorem \ref{THM-BBCC}.
\begin{cor}\label{COR-AB2}
	Let  $\widehat{m}^{\mbox{\tiny HT}}({\bf x}; \breve{\pi}_{\breve{\varphi}_n})$ be the second estimator in (\ref{CEQ5}) and suppose that (\ref{cond2}) holds.  Then, under the  conditions of Theorem \ref{THM-BBC}, for all $p\geq 2$,
	\begin{equation*} 
	E\Big|\widehat{m}^{\mbox{\tiny HT}}({\bf X}; \breve{\pi}_{\breve{\varphi}_n})- m({\bf X}) \Big|^p \,=\, \mathcal{O}\left( \sqrt{\frac{\log(\ell \vee m\vee  |\mathcal{F}_{\varepsilon_n}|)}{\ell \wedge mh^d}}\right).
	\end{equation*}
\end{cor}

\vspace{1mm}\noindent
Once again, we note that for the special case where $m=\alpha\cdot n$ and $\ell=(1-\alpha)\cdot n$, where $\alpha\in(0,1)$, under the above conditions, one finds that
\begin{equation*} 
E\Big|\widehat{m}^{\mbox{\tiny HT}}({\bf X}; \breve{\pi}_{\breve{\varphi}_n})- m({\bf X}) \Big|^p \,=\, \mathcal{O}\left( \sqrt{\frac{\log(n \vee |\mathcal{F}_{\varepsilon_n}|)}{n h^d}}\right),~~\mbox{for all $p\geq 2$}.
\end{equation*}

\section{Applications to classification with partially labeled data} \label{ACPL}
Consider the following standard two-group classification problem. Let $({\bf X}, Y)$ be a random pair, where ${\bf X}\in\mathbb{R}^d$ is a vector of covariates and $Y\in\{0, 1\}$, called the class variable or (class label), has to be predicted based on  ${\bf X}$. More specifically, the aim of classification is to find a map/function $g:\mathbb{R}\to \{0,1\}$  for which the misclassification error, i.e., 
\begin{equation} \label{Lg}
L(g) := P\{g({\bf X})\neq Y\},
\end{equation}
is as small as possible. The best classifier, also referred to as the Bayes classifier, is given by
\begin{equation}\label{Bayes}
g_{\mbox{\tiny B}}({\bf x})=\left\{
\begin{array}{ll}
1 ~ & ~ \hbox{if} ~~ m({\bf x}) \,:=\, E\big[Y|\,{\bf X}={\bf x}\big]~ > ~ \frac12 \\
0 ~ & ~ \hbox{otherwise,}
\end{array}
\right. 
\end{equation}
see, for example, Chapter 2 of Devroye et al (1996). Since, in practice,  the distribution of $({\bf X},Y)$ is virtually always unknown, finding the best classifier $g_{\mbox{\tiny B}}$ is impossible. However, suppose that we have access to $n$ iid observations (the data), $\mathbb{D}_n :=\{({\bf X}_1,Y_1), \cdots, ({\bf X}_n,Y_n)\}$, where $({\bf X}_i,Y_i)\stackrel{\mbox{\tiny iid}}{=}({\bf X},Y), ~ i=1,\cdots, n,$ and  let $\widehat{g}_{n}$ be any classifier constructed  based on the data $\mathbb{D}_n$. Also, let 
\begin{equation}
L_n(\widehat{g}_{n}) = P\big\{\widehat{g}_n({\bf X})\neq Y \big|\, \mathbb{D}_n\big\} 
\end{equation}
be the conditional misclassification error of $\widehat{g}_{n}$. Then $\widehat{g}_{n}$ is said to be weakly (strongly) Bayes consistent if $L_n(\widehat{g}_{n})\to L(g_{\mbox{\tiny B}})$ in probability (almost surely). Now, let $\widehat{m}({\bf x})$ be any estimator of the regression function $m({\bf x}) := E\big[Y|{\bf X}={\bf x}\big]$ and consider the plug-in type classifier
\begin{equation}\label{plugin}
\widehat{g}_n({\bf x})=\left\{
\begin{array}{ll}
1 ~ & ~ \hbox{if} ~~ \widehat{m}({\bf x}) ~ > ~ \frac12 \\
0 ~ & ~ \hbox{otherwise.}
\end{array}
\right. 
\end{equation}
Then, one has (see Lemma 6.1 of Devroye et al (1996))
\begin{equation}\label{Bayes-bound}
L_n(\widehat{g}_n) -L(g_{\mbox{\tiny B}}) ~\leq~ 
2 E\big[\big|\widehat{m}({\bf X}) - m({\bf X})\big|\,\big|\mathbb{D}_n\big]\,,
\end{equation}
and by the dominated convergence theorem, $E\big[L_n(\widehat{g}_n)\big] -L(g_{\mbox{\tiny B}}) \,\leq\,2 E\big|\widehat{m}({\bf X}) - m({\bf X})\big|$. Now, suppose that some of the $Y_i$'s may be missing not at random (NMAR) and consider the regression estimator
$\widehat{m}({\bf x};\widehat{\varphi}_{n})$ in (\ref{mhat5}) to be plugged in for  $m({\bf x})$ in (\ref{Bayes}). Also, denote the corresponding classifier by 
\begin{equation}\label{plugin2}
\widehat{g}_n({\bf x};\widehat{\varphi}_{n})=\left\{
\begin{array}{ll}
1 ~ & ~ \hbox{if} ~~ \widehat{m}({\bf x};\widehat{\varphi}_{n}) ~ > ~ \frac12 \\
0 ~ & ~ \hbox{otherwise.}
\end{array}
\right. 
\end{equation}
To study the asymptotic performance  of the classifier in (\ref{plugin2}), we also state the following so-called margin condition (see, for example, Audibert and Tsybakov(2007)).

\vspace{3mm}\noindent
{\bf Assumption (G)} {\it [Margin condition.]}  There exist constants $c>0$ and $\alpha>0$ such that
\begin{equation}\label{MARG}
P\left\{0\,<\, \bigg|m({\bf X})-\frac{1}{2}\bigg| \,\leq\, t  \right\} ~\leq~ c\, t^{\alpha},~~~\mbox{for all $t>0$.}
\end{equation}

\vspace{1mm}\noindent
Applications of the above margin condition in classification has been considered by many authors; see, for example,  Mammen
and Tsybakov (1999), Massart and Nédélec (2006), Audibert and Tsybakov(2007),  Tsybakov
and van de Geer (2005), Kohler and Krzy\.{z}ak (2007), and D\"{o}ring et al (2016). 

\begin{thm} \label{THM-class-1}
	Consider the classifier $\widehat{g}_n({\bf x};\widehat{\varphi}_{n})$ given by (\ref{plugin2}). If  (\ref{cond}) holds then, under the conditions of Theorem \ref{THM-B}, we have
		
	\vspace{2mm}\noindent
	\textbf{\textit{(i)}}~~~~~~~~~
	$
	P\left\{\widehat{g}_n({\bf X};\widehat{\varphi}_{n}) \neq Y\Big| \mathbb{D}_n \right\} \longrightarrow^{a.s.} P\{g_{\mbox{\tiny B}}({\bf X})\neq Y\}.
	$
	
	\vspace{2mm}\noindent
	\textbf{\textit{(ii)}}~~~~~~~~~
	$
	P\left\{\widehat{g}_n({\bf X};\widehat{\varphi}_{n}) \neq Y\right\} - P\{g_{\mbox{\tiny B}}({\bf X})\neq Y\} = \mathcal{O}\left( \left(\frac{\log(\ell\vee |\mathcal{F}_{\varepsilon_n}|)}{\ell\wedge (mh^d)}   \right)^{1/4}\right).
	$
	
	\vspace{2mm}\noindent
	\textbf{\textit{(iii)}}~ If the margin condition (\ref{MARG}) holds then
	\[
	P\left\{\widehat{g}_n({\bf X};\widehat{\varphi}_{n}) \neq Y\right\} - P\{g_{\mbox{\tiny B}}({\bf X})\neq Y\} = \mathcal{O}\left( \left(\frac{\log(\ell\vee |\mathcal{F}_{\varepsilon_n}|)}{\ell\wedge (mh^d)}   \right)^{\frac{1+\alpha}{2(2+\alpha)}}\right),
	\]
	where $\alpha$ is as in (\ref{MARG}).
\end{thm}	
Part (iii) of the above theorem shows that for large values of $\alpha$ we can obtain rates closer to $\big(\log(\ell\vee|\mathcal{F}_{\varepsilon_n}|)\big/\big[\ell\wedge (mh^d)\big]\big)^{1/2}$\, which is the same as that of the actual regression estimator (see Corollary \ref{COR-A2}).
	
\vspace{3mm}\noindent
Next, consider the Horvitz-Thompson type regression estimators given by (\ref{CEQ5}) and denote the corresponding plug-in classifiers by
\begin{equation}\label{plugin3}
\widetilde{g}^{\mbox{\tiny HT}}_n({\bf x};\widetilde{\pi})=\left\{
\begin{array}{ll}
1 ~ & ~ \hbox{if} ~~ \widehat{m}^{\mbox{\tiny HT}}({\bf x}; \widetilde{\pi}_{\widetilde{\varphi}_n}) ~ > ~ \frac12 \\
0 ~ & ~ \hbox{otherwise,}
\end{array}
\right.         ~~\mbox{and}~~~
\breve{g}^{\mbox{\tiny HT}}_n({\bf x}; \breve{\pi})=\left\{
\begin{array}{ll}
1 ~ & ~ \hbox{if} ~~ \widehat{m}^{\mbox{\tiny HT}}({\bf x}; \breve{\pi}_{\breve{\varphi}_n}) ~ > ~ \frac12 \\
0 ~ & ~ \hbox{otherwise,}
\end{array}
\right. 
\end{equation}
where $\widehat{m}^{\mbox{\tiny HT}}({\bf x}; \widetilde{\pi}_{\widetilde{\varphi}_n})$ and $\widehat{m}^{\mbox{\tiny HT}}({\bf x}; \breve{\pi}_{\breve{\varphi}_n})$ are as in (\ref{CEQ5}).
As for the asymptotic performance of the two classifiers in (\ref{plugin3}), we have the following counterpart of Theorem 
(\ref{THM-class-1}).

\begin{thm} \label{THM-class-2}
	Let $\widetilde{g}^{\mbox{\tiny HT}}_n$ and $\breve{g}^{\mbox{\tiny HT}}_n$ be the two classifiers in (\ref{plugin3}). If  (\ref{cond2}) holds then, under the conditions of Theorem \ref{THM-BBC}, we have
	
	\vspace{2mm}\noindent
	\textbf{\textit{(i)}}
	$
	P\left\{\widetilde{g}^{\mbox{\tiny HT}}_n({\bf X}; \widetilde{\pi}) \neq Y\Big| \mathbb{D}_n \right\} \rightarrow^{a.s.} P\{g_{\mbox{\tiny B}}({\bf X})\neq Y\}\,\mbox{and}~
	P\left\{\breve{g}^{\mbox{\tiny HT}}_n({\bf X}; \breve{\pi}) \neq Y\Big| \mathbb{D}_n \right\} \rightarrow^{a.s.} P\{g_{\mbox{\tiny B}}({\bf X})\neq Y\}.
	$
	
	\vspace{2mm}\noindent
	\textbf{\textit{(ii)}}~~~~~~~~~
	$
	P\left\{\breve{g}^{\mbox{\tiny HT}}_n({\bf X};\breve{\varphi}_{n}) \neq Y\right\} - P\{g_{\mbox{\tiny B}}({\bf X})\neq Y\} = \mathcal{O}\left( \left(\frac{\log(\ell\vee m \vee |\mathcal{F}_{\varepsilon_n}|)}{\ell\wedge (mh^d)}   \right)^{1/4}\right).
	$
	
	\vspace{2mm}\noindent
	\textbf{\textit{(iii)}}~ If the margin condition (\ref{MARG}) holds then
	\[
	P\left\{\breve{g}^{\mbox{\tiny HT}}_n({\bf X};\breve{\varphi}_{n}) \neq Y\right\} - P\{g_{\mbox{\tiny B}}({\bf X})\neq Y\} = \mathcal{O}\left( \left(\frac{\log(\ell\vee m\vee |\mathcal{F}_{\varepsilon_n}|)}{\ell\wedge (mh^d)}   \right)^{\frac{1+\alpha}{2(2+\alpha)}}\right),
	\]
	where $\alpha$ is as in (\ref{MARG}).
\end{thm}	
Here, we observe that for large values of $\alpha$ in part (iii) of the above theorem, one can obtain rates closer to  $\big(\log(\ell\vee m\vee |\mathcal{F}_{\varepsilon_n}|)\big/\big[\ell\wedge (mh^d)\big]\big)^{1/2}$\, which is similar to that of the  winsorized-type  regression estimator $\widehat{m}^{\mbox{\tiny HT}}({\bf x}; \breve{\pi}_{\breve{\varphi}_n})$ in (\ref{CEQ5}); see Corollary \ref{COR-AB2}.

\section{Proofs of the main results} \label{PRF}
We start by stating a number of lemmas.  In what follows, we use the notation of Section \ref{sub1} and let   $\mathcal{F}$ be a totally bounded class of functions $\varphi: [-L, L] \rightarrow (0,B],$ for some $B<\infty$. Also, for any $\varepsilon>0$, we let $\mathcal{F}_{\varepsilon}$ be any $\varepsilon$-cover of $\mathcal{F}$ (as defined in Section \ref{sub1}).  Next, for each $\varphi\in\mathcal{F}$, put
\begin{equation} \label{psieta}
\psi_k({\bf x}; \varphi) := E\left[\Delta Y^{2-k} \varphi(Y)\Big|{\bf X}={\bf x}\right]~~~\mbox{and}~~~\eta_k({\bf x}) := E\big[\Delta Y^{2-k} \big|{\bf X}={\bf x}\big],~~~\mbox{for}\,~k=1,2,
\end{equation}
and define
\begin{equation} \label{repr3}
m({\bf x;\varphi}) ~=~ \eta_1({\bf x}) + \frac{\psi_1({\bf x};\varphi)}{\psi_2({\bf x};\varphi)}\cdot \left(1- \eta_2({\bf x})\right).
\end{equation}
Also, define
\begin{eqnarray}\label{Eq1}
\widehat{L}_{m,\ell}(\varphi) &:=& \frac{1}{\ell}\, \sum_{i\in\,\boldsymbol{{\cal I}}_\ell} \frac{\Delta_i}{ \widehat{\pi}_{\varphi}({\bf X}_i, Y_i)}\,
\Big|\widehat{m}_{m}({\bf X}_i;\varphi)-Y_i\Big|^2,
\end{eqnarray}
where $\widehat{\pi}_{\varphi}({\bf z}, y)$ is as in  (\ref{FIhat1}) 
and $\widehat{m}_{m}({\bf x};\varphi)$ is given by (\ref{mhat3}). Also, put 
\begin{equation}\label{FIEP2}
\varphi_{\varepsilon} \,:=~ \argmin_{\varphi\in\mathcal{F}_{\varepsilon}}\, E\big|m({\bf X}; \varphi)-Y\big|^2 ~~~~\mbox{and} ~~~~\widehat{\varphi}_{\varepsilon} \,:=~ \argmin_{\varphi\in\mathcal{F}_{\varepsilon}}\, \widehat{L}_{m,\ell}(\varphi) .
\end{equation}

\begin{lem} \label{LEM-00}
Let $\varphi^*$ be the true (unknown) version of the function $\varphi$ in  (\ref{NonIg4}). Also, let $m({\bf x};\varphi)$ be as defined in model (\ref{repr3}). Then the regression function $m({\bf x})=E[Y|{\bf X}={\bf x}]$ can be represented as
\begin{equation}\label{repr2}
m({\bf x}) \,=\,m({\bf x};\varphi^*)\,=\,\eta_1({\bf x}) + \frac{\psi_1({\bf x};\varphi^*)}{\psi_2({\bf x};\varphi^*)}\cdot \left(1- \eta_2({\bf x})\right).
\end{equation}
where the functions $\psi_k$ and $\eta_k$, $k=1,2$, are given by (\ref{psieta}).
\end{lem}	

\vspace{2mm}\noindent
PROOF OF LEMMA \ref{LEM-00}.

\vspace{1mm}\noindent	 
The proof of this lemma is straightforward and therefore omitted. 

\hfill$\Box$

\begin{lem} \label{LEM-1} 
Let $m({\bf x};\varphi)$, $\widehat{L}_{m,\ell}(\varphi)$, $\varphi_{\varepsilon}$, and $\widehat{\varphi}_{\varepsilon}$ be as in (\ref{repr3}), (\ref{Eq1}), and (\ref{FIEP2}), respectively.  Then, under the conditions of Theorem \ref{THM-B}, we have
\begin{eqnarray}
E\left[\Big|\widehat{m}_m({\bf X}; \widehat{\varphi}_{\varepsilon}) - m({\bf X};\varphi_{\varepsilon})\Big|^2\bigg|\mathbb{D}_n\right] 
&\leq& \sup_{\varphi \in \mathcal{F}_{\varepsilon}} \,\left|E\left[\Big|\widehat{m}_m({\bf X}; \varphi) - Y\Big|^2\Big|\mathbb{D}_m\right] - \widehat{L}_{m,\ell}(\varphi) \right| \nonumber\\
&& \,+\, \sup_{\varphi \in \mathcal{F}_{\varepsilon}} \,\left| \widehat{L}_{m,\ell}(\varphi) - 
E\Big|m({\bf X}; \varphi) - Y\Big|^2 \right|
~+ ~C_1 \,\varepsilon^{1/2},  ~~~~~          \label{EQNN}
\end{eqnarray}
 where $C_1$ is a positive constant not depending on $n$ or $\varepsilon$, and $\widehat{m}_m({\bf X}; \varphi) $ is as in  (\ref{mhat3}).
\end{lem}

\vspace{2mm}\noindent
PROOF OF LEMMA \ref{LEM-1}.

\vspace{1mm}\noindent
Observe that $E\big[|\widehat{m}_m({\bf X}; \widehat{\varphi}_{\varepsilon}) - Y|^2\big|\mathbb{D}_n\big]= E\big[|\widehat{m}_m({\bf X}; \widehat{\varphi}_{\varepsilon}) - m({\bf X};\varphi_{\varepsilon})|^2\big|\mathbb{D}_n\big]+E|m({\bf X};\varphi_{\varepsilon})-Y|^2$ $+\,2 E\big[\big(\widehat{m}_m({\bf X}; \widehat{\varphi}_{\varepsilon}) - m({\bf X};\varphi_{\varepsilon})\big)\big(m({\bf X};\varphi_{\varepsilon}) - Y\big)\big|\mathbb{D}_n\big]$.
Also, let $\varphi^*$ be as in (\ref{PHIstar}) and note that
\begin{eqnarray*}
&& E\left[\Big(\widehat{m}_m({\bf X}; \widehat{\varphi}_{\varepsilon}) - m({\bf X};\varphi_{\varepsilon})\Big)\Big(m({\bf X};\varphi_{\varepsilon}) - Y\Big)\bigg|\mathbb{D}_n\right] \\
&& ~= E\left[\Big(\widehat{m}_m({\bf X}; \widehat{\varphi}_{\varepsilon}) - m({\bf X};\varphi_{\varepsilon})\Big)\Big(m({\bf X};\varphi_{\varepsilon}) - 
m({\bf X};\varphi^*) + m({\bf X};\varphi^*) - Y\Big)\bigg|\mathbb{D}_n\right]\\
&& ~= E\left[\Big(\widehat{m}_m({\bf X}; \widehat{\varphi}_{\varepsilon}) - m({\bf X};\varphi_{\varepsilon})\Big)\Big(m({\bf X};\varphi_{\varepsilon}) - 
m({\bf X};\varphi^*) \Big)\bigg|\mathbb{D}_n\right],  
\end{eqnarray*}
where we have used the fact that in view of (\ref{repr2}), $E[Y|{\bf X}={\bf x}]:=m({\bf x})=m({\bf x}; \varphi^*)$. Therefore
\begin{eqnarray}
&&
E\left[\Big|\widehat{m}_m({\bf X}; \widehat{\varphi}_{\varepsilon}) - m({\bf X};\varphi_{\varepsilon})\Big|^2\bigg|\mathbb{D}_n\right] \nonumber\\
&& ~=~\left\{E\left[\Big|\widehat{m}_m({\bf X}; \widehat{\varphi}_{\varepsilon}) - Y\Big|^2\bigg|\mathbb{D}_n\right] - E\Big|m({\bf X}; \varphi_{\varepsilon}) -Y\Big|^2\right\}\nonumber\\
&&~~~~-\,2E\left[\Big(\widehat{m}_m({\bf X}; \widehat{\varphi}_{\varepsilon}) - m({\bf X};\varphi_{\varepsilon})\Big)\Big(m({\bf X};\varphi_{\varepsilon}) - 
m({\bf X};\varphi^*) \Big)\bigg|\mathbb{D}_n\right]\nonumber\\
&&~:=\, {\bf I}_n +{\bf I\!I}_n. \label{EQ4}
\end{eqnarray}
Now, observe that
\begin{eqnarray*}
{\bf I}_n &=& E\left[\Big|\widehat{m}_m({\bf X}; \widehat{\varphi}_{\varepsilon}) - Y\Big|^2\bigg|\mathbb{D}_n\right] - 
\inf_{\varphi\in\mathcal{F}_{\varepsilon}}\, E\big|m({\bf X}; \varphi)-Y\big|^2\\
&=& \sup_{\varphi\in\mathcal{F}_{\varepsilon}}\, \bigg\{E\left[\Big|\widehat{m}_m({\bf X}; \widehat{\varphi}_{\varepsilon}) - Y\Big|^2\bigg|\mathbb{D}_n\right] -   
\widehat{L}_{m,\ell}(\varphi) + \widehat{L}_{m,\ell}(\varphi) - \widehat{L}_{m,\ell}(\widehat{\varphi}_{\varepsilon}) \\
&&~~+\, \widehat{L}_{m,\ell}(\widehat{\varphi}_{\varepsilon})  - E\big|m({\bf X}; \varphi)-Y\big|^2 \bigg\},~~~~\mbox{(where $ \widehat{L}_{m,\ell}(\varphi)$ is as in (\ref{Eq1}))}\\
&\leq& \bigg(E\left[\Big|\widehat{m}_m({\bf X}; \widehat{\varphi}_{\varepsilon}) - Y\Big|^2\bigg|\mathbb{D}_n\right]  -  \widehat{L}_{m,\ell}(\widehat{\varphi}_{\varepsilon}) \bigg) + \sup_{\varphi\in\mathcal{F}_{\varepsilon}}\, \left|\widehat{L}_{m,\ell}(\varphi) - E\big|m({\bf X}; \varphi)-Y\big|^2 \right|,
\end{eqnarray*}
where the last line follows since $\widehat{L}_{m,\ell}(\widehat{\varphi}_{\varepsilon})\leq \widehat{L}_{m,\ell}(\varphi)$ holds for all $\varphi\in \mathcal{F}_{\varepsilon}$ (because of the definition of $\widehat{\varphi}_{\varepsilon}$ in  (\ref{FIEP2})). Therefore, 
\begin{equation}
\big|{\bf I}_n\big| \,\leq\, \sup_{\varphi\in\mathcal{F}_{\varepsilon}}\, \Bigg|E\left[\Big|\widehat{m}_m({\bf X}; \varphi) - Y\Big|^2\bigg|\mathbb{D}_m\right]  -  \widehat{L}_{m,\ell}(\varphi) \Bigg| + \sup_{\varphi\in\mathcal{F}_{\varepsilon}}\, \Bigg|\widehat{L}_{m,\ell}(\varphi) - E\big|m({\bf X}; \varphi)-Y\big|^2\Bigg|,~\label{EQ4Q}
\end{equation}
where the conditioning on $\mathbb{D}_m$ in the above expression  reflects the fact that $\widehat{m}_m({\bf X};\varphi)$ depends on $\mathbb{D}_m$ only (and not the entire data $\mathbb{D}_n$).  Furthermore, the term ${\bf I\!I}_n$ in (\ref{EQ4}) can be bounded as follows. 
\begin{eqnarray}
\big|{\bf I\!I}_n\big| &\leq& 2E\left[\Big|\widehat{m}_m({\bf X}; \widehat{\varphi}_{\varepsilon}) - m({\bf X};\varphi_{\varepsilon})\Big|\cdot\Big|
m({\bf X};\varphi_{\varepsilon}) - m({\bf X};\varphi^*) \Big|\bigg|\mathbb{D}_n\right]\nonumber\\
&\leq& 6L\cdot E\Big| m({\bf X};\varphi_{\varepsilon}) - m({\bf X};\varphi^*) \Big| ~\leq~ 6L\,\sqrt{E\big| m({\bf X};\varphi_{\varepsilon}) - m({\bf X};\varphi^*) \big|^2}.~~ 
\label{EQ4C}
\end{eqnarray}
But, using the identity $E\big| m({\bf X};\varphi_{\varepsilon}) - Y\big|^2=E\big| m({\bf X};\varphi^*) - Y\big|^2 + E\big| m({\bf X};\varphi_{\varepsilon}) - m({\bf X};\varphi^*) \big|^2$, we have 
\begin{eqnarray}
E\big| m({\bf X};\varphi_{\varepsilon}) - m({\bf X};\varphi^*) \big|^2 &=& \inf_{\varphi\in\mathcal{F}_{\varepsilon}} E\big| m({\bf X};\varphi) - Y\big|^2
- E\big| m({\bf X};\varphi^*) - Y\big|^2\nonumber\\
&=& \inf_{\varphi\in\mathcal{F}_{\varepsilon}} E\big| m({\bf X};\varphi) - m({\bf X};\varphi^*) \big|^2 
\nonumber\\   &\leq&   2L \inf_{\varphi\in\mathcal{F}_{\varepsilon}} E\big| m({\bf X};\varphi) - m({\bf X};\varphi^*) \big|.  \label{EQ4B}
\end{eqnarray}
Now let $\varphi^{\dagger}\in \mathcal{F}_{\varepsilon}$ be such that $\varphi^*\in B(\varphi^{\dagger}, \varepsilon)$; such a $\varphi^{\dagger}\in \mathcal{F}_{\varepsilon}$ exists because  $\varphi^*\in \mathcal{F}$ and $\mathcal{F}_{\varepsilon}$ is an $\varepsilon$-cover of $\mathcal{F}$.  Then, in view of Lemma \ref{LEM-2} and the fact that the  right side of (\ref{EQ4B}) is an infimum, one finds
\begin{eqnarray}
\mbox{(Rght side of (\ref{EQ4B})})&\leq& 2L\cdot  E\big| m({\bf X};\varphi^{\dagger}) - m({\bf X};\varphi^*) \big|  \,\leq\, 2LC \sup_{-L\leq y\leq L} \big|\varphi^{\dagger}(y) - \varphi^*(y) \big| \nonumber \\
&\leq& 2LC \cdot \varepsilon ~~~~\big(\mbox{because $\varphi^*\in B(\varphi^{\dagger}, \varepsilon)$}\big),   \label{EEQ2}
\end{eqnarray}
where $C$ is as in Lemma \ref{LEM-2}. Therefore, by (\ref{EQ4C}) and (\ref{EQ4B}), we have
\begin{equation}
\big|{\bf I\!I}_n\big|\, \leq~  6L\,\sqrt{2LC\cdot \varepsilon} \,~=:~ C_1\sqrt{\varepsilon\,}\,.  \label{EQ4D}
\end{equation}
Now Lemma \ref{LEM-1} follows from (\ref{EQ4}), (\ref{EQ4Q}), and (\ref{EQ4D}).

\hfill $\Box$

\begin{lem}\label{LEM-2}
Let $m({\bf x};\varphi_j)$, $j=1,2$,  be as defined in (\ref{repr3}), where $\varphi_j:\, [-L,L] \to (0, B]$ for some positive number $B$. Then, under assumption (A4), one has
\begin{equation*}
E\Big|m({\bf X};\varphi_1) - m({\bf X};\varphi_2)\Big| ~\leq~  C\cdot \sup_{-L\leq y \,\leq L}\big|\varphi_1(y)-\varphi_2(y)\big|,
\end{equation*}
where  the constant $C>0$ can be taken to be $C=2L/\varrho_0$, with $\varrho_0$ as in assumption (A4).
\end{lem}

\vspace{2mm}\noindent
PROOF OF LEMMA \ref{LEM-2}.

\vspace{1mm}\noindent
Let
$
S_j({\bf x})=E[\Delta Y\,\varphi_j(Y)|{\bf X}= {\bf x}]~\mbox{and}~ T_j({\bf x})=E[\Delta\,\varphi_j(Y)|{\bf X}= {\bf x}],~j=1,2,
$
and observe that 
\begin{eqnarray*}
\Big| m({\bf x};\varphi_1) - m({\bf x};\varphi_2) \Big| 
&=&\left| \frac{- S_1({\bf x})}{T_1({\bf x})}\cdot 
\frac{T_1({\bf x})-T_2({\bf x})}{T_2({\bf x})}
+ \frac{S_1({\bf x})-S_2({\bf x})}{T_2({\bf x})} \right|\cdot E\big[1-\delta|{\bf X}={\bf x}\big]\nonumber\\
&\leq&\frac{1}{T_2({\bf x})}\big\{L\big|T_1({\bf x})-T_2({\bf x})\big| +\big|S_1({\bf x})-S_2({\bf x})\big|
\big\}. 
\end{eqnarray*}
But,  $|S_1({\bf x}) - S_2({\bf x})| \leq\, E\left[|\delta\,Y|\cdot\big|\varphi_1(Y)-\varphi_2(Y)\big|\, \big|{\bf X}={\bf x} \right] \leq\, L\,\sup_{-L\leq y \,\leq L}\big|\varphi_1(y)-\varphi_2(y)\big|$. Similarly, $|T_1({\bf x})-T_2(\boldsymbol{x})|\,\leq\, \sup_{-L\leq y \,\leq L}\big|\varphi_1(y)-\varphi_2(y)\big|$.
On the other hand, by Assumption (D), we have $T_2({\bf x})\geq\varrho_{0}>0$,  for $\mu$--a.e.\,${\bf x}$. Therefore
\[
\Big| m({\bf x};\varphi_1) - m({\bf x};\varphi_2) \Big|  \,\leq \, (2L/\varrho_0) \,\sup_{-L\leq y \,\leq L}\big|\varphi_1(y)-\varphi_2(y)\big|
\]
The lemma follows now by integrating both sides of this inequality with respect to $\mu({d \bf x})$.

\hfill $\Box$

\begin{lem}\label{L-2}
	Let $\mathcal{K}$ be a regular kernel. Also, let  $\mu$ be any probability measure on the Borel sets of $\mathbb{R}^d$. Then there is a positive constant $\rho(\mathcal{K})$, depending on the kernel $\mathcal{K}$ but not $n$, such that for every $h>0$,  
	$$
	\sup_{{\bf u}\in\mathbb{R}^d}\,\int\frac{ \mathcal{K}(({\bf x}-{\bf u})/h)}{E[\mathcal{K}(({\bf x} -{\bf X})/h)]}\,\mu(d{\bf x})\,\leq \,\rho(\mathcal{K})\,.$$
\end{lem}
\noindent
{\bf Proof.} The proof of this lemma appears in Devroye and  Krzy\`{z}ak (1989; Lemma 1).  

$\hfill\square$

\begin{lem}\label{L-3}
	Let $({\bf U},Y), ({\bf U}_1,Y_1),\dots, ({\bf U}_n,Y_n)$ be iid $\mathbb{R}^d\times [-A,A]$-valued random vectors for some $A\in[0,\infty)$. Also, let\, $m({\bf u})= E[Y|{\bf U}$=${\bf u}]$ be the regression function and define  quantity $\widetilde{m}_n({\bf u})=\sum_{i=1}^n Y_i\, K\big(({\bf u}-{\bf U}_i)/h\big)\big/
	\big\{nE[\mathcal{K}(({\bf u}-{\bf U})/h)]\big\}$, where $\mathcal{K}$ is a regular kernel.  If $h\to 0$ and $n h^d\to \infty$, as $n\to\infty$, then for every $t>0$ and large enough $n$, 
	$$
	P\left\{\displaystyle\int\Big|\widetilde{m}_n({\bf u})
	- m({\bf u})\Big|\,\mu(d{\bf u})>t\right\}\,\leq\, e^{-n t^2/(64 A^2\,  \rho^2(\mathcal{K}))}\,
	$$ 
	where   $\mu$ is the probability measure of ${\bf U}$, and $\rho(\mathcal{K})$ is as in Lemma \ref{L-2}.
\end{lem}

\vspace{1mm}\noindent
{\bf Proof.} The proof can be found in Gy\"{o}rfi et al. (2002; Sec. 23).

\hfill $\Box$

\vspace{0.2mm}\noindent
PROOF OF THEOREM \ref{THM-A}

\vspace{1mm}\noindent
Firt observe that for every $p\geq$\,1,  
$| \widehat{m}_{n,\widehat{\gamma}}({\bf x}) - m({\bf x})|^p \leq \{| \widehat{m}_{n,\widehat{\gamma}}({\bf x})|+| m({\bf x})|\}^{p-1}| \widehat{m}_{n,\widehat{\gamma}}({\bf x})   - m({\bf x})|  \leq  (3L)^{p-1}| \widehat{m}_{n,\widehat{\gamma}}({\bf x})   - m({\bf x})|$. 
Therefore, we only need to prove the theorem for the case of $p=$\,1. 
To this end, let $\eta_k({\bf x}, t)$ and $\widehat{\eta}_k({\bf x}, t)$, $k=1,2$, $t\in\mathbb{R}$, be the quantities defined in (\ref{ETA12}) and (\ref{ETA12.hat}), respectively.
Then it is straightforward to show that in view of  (\ref{mhat}) and  (\ref{repr}), and the fact that $|\widehat{\eta}_1({\bf x}, \widehat{\gamma})/\widehat{\eta}_2({\bf x}, \widehat{\gamma})|\leq L$, one has
\begin{eqnarray}
\big| \widehat{m}_{n,\widehat{\gamma}}({\bf x}) - m({\bf x})\big| &\leq& \big|\widehat{\eta}_1({\bf x}, 0) - \eta_1({\bf x}, 0)\big|
+ \left|\frac{\widehat{\eta}_1({\bf x}, \widehat{\gamma})}{\widehat{\eta}_2({\bf x}, \widehat{\gamma})} - \frac{\eta_1({\bf x}, \gamma^*)}{\eta_2({\bf x}, \gamma^*)}\right| + L\cdot\big|\widehat{\eta}_2({\bf x}, 0) - \eta_2({\bf x}, 0)\big|.~~~~~ \label{EQQ1}
\end{eqnarray}
But the first and third terms on the right side of (\ref{EQQ1}) can be immediately bounded using the classical result of Devroye and Krzy\`{z}ak (1989). More specifically, for every $t>0$ and $n$ large enough,
\begin{equation}\label{EQ2Bound}
P\left\{\int\big|\widehat{\eta}_1({\bf x}, 0) - \eta_1({\bf x}, 0)\big|\mu(d{\bf x}) >t \right\}\leq e^{-a_1 n}~\,\mbox{and}~
P\left\{\int L\big|\widehat{\eta}_2({\bf x}, 0) - \eta_2({\bf x}, 0)\big|\mu(d{\bf x}) >t \right\}\leq e^{-a_2 n},
\end{equation}
where $a_1$ and $a_2$ are positive constants that depend on $t$ but not $n$. To deal with the middle term on the right side of (\ref{EQQ1}), we note that it can be written as 
\begin{eqnarray*}
&& \left|\frac{1}{\eta_2({\bf x}, \gamma^*)} \left[\frac{\widehat{\eta}_1({\bf x}, \widehat{\gamma})}{\widehat{\eta}_2({\bf x}, \widehat{\gamma})} \cdot
\big(\eta_2({\bf x}, \gamma^*) - \widehat{\eta}_2({\bf x}, \widehat{\gamma}) \big)+ \big(\widehat{\eta}_1({\bf x}, \widehat{\gamma}) - \eta_1({\bf x}, \gamma^*) \big) \right]\right| \nonumber\\
&&  \leq\, \frac{1}{\pi_{\mbox{\tiny min}}\exp\{-L|\gamma^*|\}}
\Big[ \big|\widehat{\eta}_1({\bf x}, \widehat{\gamma}) - \widehat{\eta}_1({\bf x}, \gamma^*)\big| +  \big|\widehat{\eta}_1({\bf x}, \gamma^*) - \eta_1({\bf x}, \gamma^*)\big| + L\big|\widehat{\eta}_2({\bf x}, \widehat{\gamma}) - \widehat{\eta}_2({\bf x}, \gamma^*)\big| \nonumber\\
&& ~~~~~ + L \big|\widehat{\eta}_2({\bf x}, \gamma^*) - \eta_2({\bf x}, \gamma^*)\big| \Big],
\end{eqnarray*}
where the above inequality follows from Assumption (A) with the simple fact that $\eta_2({\bf X}, \gamma^*)$ $=E\left(E\left[\Delta \exp\{\gamma^* Y\}\big| {\bf X}, Y\right]\big|{\bf X}\right)=E\left[\exp\{\gamma^* Y\} \pi_{\gamma^*}({\bf X},Y)\big| {\bf X}\right]\,\geq \,\inf_{{\bf z},y}\, \pi_{\gamma^*}({\bf z}, y) \exp(-|\gamma^*|L)$, together with the observation that $\widehat{\eta}_1({\bf x}, \widehat{\gamma})\big/
\widehat{\eta}_2({\bf x}, \widehat{\gamma}) \leq L$. Consequently, for every $t>0$, the integral of the middle term on the right side of (\ref{EQQ1}) can be dealt with as follows
\begin{eqnarray}
P\left\{\int   \left|\frac{\widehat{\eta}_1({\bf x}, \widehat{\gamma})}{\widehat{\eta}_2({\bf x}, \widehat{\gamma})} - \frac{\eta_1({\bf x}, \gamma^*)}{\eta_2({\bf x}, \gamma^*)}\right|  \mu(d{\bf x}) >t \right\} 
&\leq& P \left\{\int    \big|\widehat{\eta}_1({\bf x}, \widehat{\gamma}) - \widehat{\eta}_1({\bf x}, \gamma^*)\big|  \mu(d{\bf x}) >\frac{\pi_{\mbox{\tiny min}} t}{4\exp\{L|\gamma^*|\}} \right\}\nonumber\\
&& + P \left\{\int    \big|\widehat{\eta}_1({\bf x}, \gamma^*) - \eta_1({\bf x}, \gamma^*)\big|  \mu(d{\bf x}) >\frac{\pi_{\mbox{\tiny min}} t}{4\exp\{L|\gamma^*|\}} \right\}\nonumber\\
&& + P \left\{\int    \big|\widehat{\eta}_2({\bf x}, \widehat{\gamma}) - \widehat{\eta}_2({\bf x}, \gamma^*)\big|  \mu(d{\bf x}) >\frac{\pi_{\mbox{\tiny min}} t}{4L\exp\{L|\gamma^*|\}} \right\}\nonumber\\
&& + P \left\{\int    \big|\widehat{\eta}_2({\bf x}, \gamma^*) - \eta_2({\bf x}, \gamma^*)\big|  \mu(d{\bf x}) >\frac{\pi_{\mbox{\tiny min}} t}{4L\exp\{L|\gamma^*|\}} \right\}\nonumber\\
&& := \sum_{k=1}^4 \mathcal{P}_{nk}(t).  \label{EQQ2}
\end{eqnarray}
To deal with the first term in (\ref{EQQ2}), i.e., the term $\mathcal{P}_{n1}(t)$, put
\begin{eqnarray*}
\Gamma'_n({\bf x}) &=& \sum_{i=1}^n\Delta_i Y_i \,\big(\exp\left\{\widehat{\gamma} \,Y_i\right\} -\exp\left\{\gamma^* \,Y_i\right\} \big)\mathcal{K}(({\bf x} -{\bf X}_i)/h) \Big/ nE\big[\mathcal{K}(({\bf x}-{\bf X})/h)\big]\\
\Gamma''_n({\bf x}) &=& \frac{\sum_{i=1}^n\Delta_i Y_i \,\big(\exp\left\{\widehat{\gamma} \,Y_i\right\} -\exp\left\{\gamma^* \,Y_i\right\} \big)\mathcal{K}(({\bf x} -{\bf X}_i)/h)}{\sum_{i=1}^n\mathcal{K}(({\bf x}-{\bf X}_i)/h)}\cdot \left[ \frac{\sum_{i=1}^n\mathcal{K}(({\bf x}-{\bf X}_i)/h)}{nE\big[\mathcal{K}(({\bf x}-{\bf X})/h)\big]} -1 \right]
\end{eqnarray*}

\vspace{1mm}\noindent
and observe that 
\begin{eqnarray} \label{Pn1B}
\mathcal{P}_{n1}(t)  &\leq&   P \left\{\int   \big|\Gamma'_n({\bf x}) \big|    \mu(d{\bf x}) >\frac{\pi_{\mbox{\tiny min}} t}{8\exp\{L|\gamma^*|\}} \right\}
+  P \left\{\int   \big|\Gamma''_n({\bf x}) \big|    \mu(d{\bf x}) >\frac{\pi_{\mbox{\tiny min}} t}{8\exp\{L|\gamma^*|\}} \right\}.~~~
\end{eqnarray}
Furthermore,
\begin{eqnarray}
\int   \big|\Gamma'_n({\bf x}) \big|  \mu(d{\bf x}) 
&\leq& 
\sup_{\bf z} \int \frac{\mathcal{K}(({\bf x}-{\bf z})/h)}{E\big[\mathcal{K}(({\bf x}-{\bf X})/h)\big]} \mu(d{\bf x})  \cdot \frac{1}{n}\,\sum_{i=1}^n \Big|\Delta_i Y_i \,\Big(\exp\left\{\widehat{\gamma} \,Y_i\right\} -\exp\left\{\gamma^* \,Y_i\right\} \Big)\Big| \nonumber\\
&\leq&
\frac{L\rho(\mathcal{K})}{n} \sum_{i=1}^n \Big|\exp\left\{\widehat{\gamma} \,Y_i\right\} -\exp\left\{\gamma^* \,Y_i\right\}\Big|,~~~\mbox{(by Lemma \ref{L-2})} \nonumber\\
&\leq& 
n^{-1}L\,\rho(\mathcal{K}) \sum_{i=1}^n \Big|\big(\widehat{\gamma}-\gamma^*\big) Y_i\,\exp\big\{\overline{\gamma} Y_i  \big\} \Big|,~~~\mbox{(via a one-term Taylor expansion)},   \nonumber\\
&& ~~\mbox{(where $\overline{\gamma}$ is a point in the interior of the line segment joining $\widehat{\gamma}$ and $\gamma^*$),}\nonumber\\
&\leq& 
n^{-1}L^2\,\rho(\mathcal{K}) \big|\widehat{\gamma}-\gamma^*\big|\cdot\sum_{i=1}^n \exp\Big\{\big|\overline{\gamma}-\gamma^*\big|L+ \big|\gamma^*|Y_i  \Big\} \label{LAB1}
\end{eqnarray}
Therefore, using the fact that $|\overline{\gamma}-\gamma^*| \leq |\widehat{\gamma}-\gamma^*|$, one finds, for every constants $t>0$ and $C_0>0$,
\begin{eqnarray}
&& 
P\left\{ \int   \big|\Gamma'_n({\bf x}) \big|  \mu(d{\bf x})  > \frac{\pi_{\mbox{\tiny min}} t}{8\exp\{L|\gamma^*|\}}\right\}\nonumber\\
&&
\leq P\left\{ \big|\widehat{\gamma}-\gamma^*\big| \exp\big\{ \big|\widehat{\gamma}-\gamma^*\big| L \big\}\cdot \frac{1}{n} \sum_{i=1}^n \exp\big\{ \big|\gamma^* \big| Y_i\big\}   >  \frac{\pi_{\mbox{\tiny min}} t}{8 L^2 \,\rho(\mathcal{K})  \exp\{L|\gamma^*|\}}\right\}\nonumber\\
&&
\leq P\Bigg\{ \left[\big|\widehat{\gamma}-\gamma^*\big| \exp\big\{ \big|\widehat{\gamma}-\gamma^*\big| L \big\}\cdot \frac{1}{n} \sum_{i=1}^n \exp\big\{ \big|\gamma^* \big| Y_i\big\}  >\,  \frac{\pi_{\mbox{\tiny min}} t}{8 L^2 \rho(\mathcal{K})  \exp\{L|\gamma^*|\}}\right]\nonumber\\
&&~~~~~~~\cap\Big\{|\widehat{\gamma}-\gamma^*|\leq C_0 \Big\}\Bigg\}\,+\, P\Big\{|\widehat{\gamma}-\gamma^*|>C_0 \Big\}\nonumber\\[2pt]
&&
\leq n\cdot P\left\{ \exp\big\{ \big|\gamma^* \big| Y_1\big\}  \,>\,  \frac{\pi_{\mbox{\tiny min}} t}{8 L^2 C_0\,\rho(\mathcal{K})  \exp\big\{(|\gamma^*|+C_0)L\big\}}\right\} + P\Big\{|\widehat{\gamma}-\gamma^*|>C_0 \Big\} \label{nPP1}\\[4pt]
&& = 0+P\Big\{|\widehat{\gamma}-\gamma^*|>C_0 \Big\}, ~\mbox{\big(for any $C_0$ satisfying $4 L^2 C_0\,\rho(\mathcal{K})  \exp\big\{(|\gamma^*|+C_0)L\big\} < \pi_{\mbox{\tiny min}} t$\big)}~~~~ \label{nP2}
\end{eqnarray}
where the last line follows because the random variable $\exp\{ |\gamma^* | Y_1\}$ is bounded by $\exp\{|\gamma^*|L\}$, which implies that taking $C_0$ small enough forces the first probability statement in (\ref{nP1})  to be zero. As for the term $\Gamma''_n({\bf x})$, we note that in view of (\ref{LAB1}) and  the observation that   $|\overline{\gamma}-\gamma^*| \leq |\widehat{\gamma}-\gamma^*|$, one obtains
\begin{equation*}
\int   \big|\Gamma''_n({\bf x}) \big|  \mu(d{\bf x})   \,\leq\,    L^2 \big|\widehat{\gamma}-\gamma^*\big| \, \max_{1\leq i \leq n} \exp\Big\{\big|\widehat{\gamma}-\gamma^*\big|L+ \big|\gamma^*|Y_i  \Big\}\cdot \int \bigg|\frac{\sum_{j=1}^n K(({\bf x} -{\bf X}_j)/h)} {nE\big[K(({\bf x} -{\bf X})/h)\big]}-1\bigg|\,\mu(d{\bf x}).
\end{equation*}
Now, observe that
\begin{eqnarray}
&& 
P\left\{ \int   \big|\Gamma''_n({\bf x}) \big|  \mu(d{\bf x}) \,>\, \frac{\pi_{\mbox{\tiny min}} t}{8\exp\{L|\gamma^*|\}}\right\} \nonumber\\
&& \leq 
P\left\{ \big|\widehat{\gamma}-\gamma^*\big| \, \exp\big\{\big|\widehat{\gamma}-\gamma^*\big|L  \big\} \cdot
\int   \bigg|\frac{\sum_{j=1}^n K(({\bf x} -{\bf X}_j)/h)} {nE\big[K(({\bf x} -{\bf X})/h)\big]}-1\bigg| \, \mu(d{\bf x}) \,>\, \frac{\pi_{\mbox{\tiny min}} t}{8 L^2 \exp\{2L|\gamma^*|\}}\right\}\nonumber\\
&& \leq 
P\left\{  \int   \bigg|\frac{\sum_{j=1}^n K(({\bf x} -{\bf X}_j)/h)} {nE\big[K(({\bf x} -{\bf X})/h)\big]}-1\bigg|  \mu(d{\bf x})   >  \frac{\pi_{\mbox{\tiny min}} t}{8 L^2 C_0  \exp\big\{(2|\gamma^*|\mbox{$+$}C_0)L\big\}}\right\} + P\big\{|\widehat{\gamma}-\gamma^*|>C_0 \big\} \label{nP1}\nonumber\\[4pt]
&& \leq \exp\left\{\frac{-n\, \pi^2_{\mbox{\tiny min}} t^2}{64^2 L^4 C_0^2 \rho^2(\mathcal{K}) \cdot \exp\big\{2\big(2|\gamma^*| +C_0\big)L  \big\}}
\right\} + P\Big\{|\widehat{\gamma}-\gamma^*|>C_0 \Big\},\label{nP3}
\end{eqnarray}
for large $n$, by Lemma \ref{L-3}, where $C_0$ is  as in (\ref{nPP1}); here, we have used Lemma \ref{L-3} with $m({\bf u})=1$ and $Y_i=1$ for all $i=1,\cdots, n.$ Putting together (\ref{Pn1B}), (\ref{nP2}),  and  (\ref{nP3}), we find 
\begin{equation} 
\mathcal{P}_{n1}(t)  \,\leq \, \exp\{-C_2 n ^2 t^2\} + 2 P\Big\{|\widehat{\gamma}-\gamma^*|>C_0 \Big\},      \label{nPP4}
\end{equation}
for $n$ large enough, where $C_2$ is a positive constant not depending on $n$. It is a simple exercise to show that the term  $\mathcal{P}_{n3}(t)$ in   (\ref{EQQ2}) can also be bounded by the right side of (\ref{nPP4}).  Furthermore, as in (\ref{EQ2Bound}), once again we can invoke the result of Devroye and Krzy\`{z}ak (1989) to conclude that $\mathcal{P}_{n2}(t)  \leq e^{-a_3 n}$ and $\mathcal{P}_{n4}(t)  \leq e^{-a_4 n}$, for $n$ large enough, where $a_3$ and $a_4$ are positive constants not depending on $n$.  These observations in conjunction with (\ref{nPP4}), (\ref{EQQ2}), (\ref{EQ2Bound}), and (\ref{EQQ1}) complete the proof of Theorem  \ref{THM-A}.

\hfill $\Box$

\vspace{3mm}\noindent
PROOF OF THEOREM \ref{THM-B}

\vspace{1.5mm}\noindent
{\it Part (i)}

\vspace{1mm}\noindent
To proceed with the proof, first note that for each $i\in \boldsymbol{{\cal I}}_\ell\,,$ we have
\begin{eqnarray}
\frac{\Delta_i}{\widehat{\pi}_{\varphi}({\bf Z}_i, Y_i)} \Big|\widehat{m}_m({\bf X}_i; \varphi) -Y_i\Big|^2 &=&
 \frac{\Delta_i\,\big|\widehat{m}_m({\bf X}_i; \varphi) -Y_i\big|^2}{\pi_{\varphi}({\bf Z}_i, Y_i)} \nonumber \\
&& - \, \Delta_i\,\big|\widehat{m}_m({\bf X}_i; \varphi) -Y_i\big|^2
\left[\frac{1}{\pi_{\varphi}({\bf Z}_i, Y_i)}- \frac{1}{\widehat{\pi}_{\varphi}({\bf Z}_i, Y_i)}\right]. ~~~~~\label{STAR1}
\end{eqnarray}
Therefore, in view of (\ref{STAR1}) and the definition of $\widehat{L}_{m,\ell}(\varphi)$ in (\ref{Eq1}), one finds  for every $\beta>0$ 
\begin{eqnarray}
&& P\left\{ \sup_{\varphi \in \mathcal{F}_{\varepsilon}} \,\left| \widehat{L}_{m,\ell}(\varphi) - E\left[\Big|\widehat{m}_m({\bf X}; \varphi) - Y\Big|^2\Big|\mathbb{D}_m\right] \right|  > \beta \right\} \nonumber\\
&&~\leq\,   P\left\{ \sup_{\varphi \in \mathcal{F}_{\varepsilon}} \,\left|\ell^{-1} \sum_{i\in\,\boldsymbol{{\cal I}}_\ell} \frac{\Delta_i \,
\big|\widehat{m}_{m}({\bf X}_i;\varphi)-Y_i\big|^2}{\pi_{\varphi}({\bf Z}_i, Y_i)}  - E\left[\Big|\widehat{m}_m({\bf X}; \varphi) - Y\Big|^2\Big|\mathbb{D}_m\right] \right| > \frac{\beta}{2}\right\} \nonumber\\
&&~~+\, P\left\{ \sup_{\varphi \in \mathcal{F}_{\varepsilon}} \,\left|\ell^{-1} \sum_{i\in\,\boldsymbol{{\cal I}}_\ell} \Delta_i \,
\big|\widehat{m}_{m}({\bf X}_i;\varphi)-Y_i\big|^2 \left[\frac{1}{\pi_{\varphi}({\bf Z}_i, Y_i)}-\frac{1}{\widehat{\pi}_{\varphi}({\bf Z}_i, Y_i)}\right]\right|>\frac{\beta}{2}\right\} \nonumber\\[4pt]
&& ~:=\, S_n(1)+S_n(2).  \label{EQ5}
\end{eqnarray}
But for each $i\in \boldsymbol{{\cal I}}_\ell\,,$
\begin{eqnarray*}
E\left[  \frac{\Delta_i}{\pi_{\varphi}({\bf Z}_i, Y_i)}\,\Big|\widehat{m}_m({\bf X}_i; \varphi) -Y_i\Big|^2\,\Bigg| \mathbb{D}_m\right]
&=&E\left[  \frac{\Big|\widehat{m}_m({\bf X}_i; \varphi) -Y_i\Big|^2}{\pi_{\varphi}({\bf Z}_i, Y_i)}\, E\Big(\Delta_i\Big| \mathbb{D}_m,{\bf X}_i, Y_i\Big)
\,\Bigg| \mathbb{D}_m\right]\nonumber\\
&=& E\left[\Big|\widehat{m}_m({\bf X}; \varphi) -Y\Big|^2\,\Big| \mathbb{D}_m\right]. 
\end{eqnarray*}
Furthermore, conditional on the training set $\mathbb{D}_m$, the terms $\Delta_i \,
\big|\widehat{m}_{m}({\bf X}_i;\varphi)-Y_i\big|^2
\big/\pi_{\varphi}({\bf Z}_i, Y_i)$, $i\in\,\boldsymbol{{\cal I}}_\ell$, are independent bounded random variables, taking values in $\big[0,\, (3L)^2/\pi_{\mbox{\tiny min}}\big]$. Therefore,
\begin{eqnarray}
 S_n(1) &\leq & \big|\mathcal{F}_{\varepsilon}\big|\sup_{\varphi \in \mathcal{F}_{\varepsilon}} \, E \Bigg[ P\Bigg\{ \Bigg|\ell^{-1} \sum_{i\in\,\boldsymbol{{\cal I}}_\ell} \frac{\Delta_i \,
\big|\widehat{m}_{m}({\bf X}_i;\varphi)-Y_i\big|^2}{\pi_{\varphi}({\bf Z}_i, Y_i)} \nonumber\\
&& ~~~~~~~~~~~~~~~~~~~~~~~~ -\, E\left[\Big|\widehat{m}_m({\bf X}; \varphi) - Y\Big|^2\Big|\mathbb{D}_m\right] \Bigg| > \frac{\beta}{2}\,\Bigg| \mathbb{D}_m\Bigg\}\Bigg] \nonumber\\
&\leq& 2\,\big|\mathcal{F}_{\varepsilon}\big|\,  \exp\big\{-\pi^2_{\mbox{\tiny min}} \ell \beta^2/(162 L^4)\big\},~~~(\mbox{via Hoeffding's inequality}). \label{EQ6}
\end{eqnarray}
To deal with the term $S_n(2)$ in (\ref{EQ5}),  let $\widetilde{\psi}_{m}({\bf Z}_i;\varphi)$ and $\widetilde{\eta}_{m}({\bf Z}_i)$ be as in (\ref{ETPS}) and observe that in view of (\ref{NonIg4}), (\ref{expgz2}), (\ref{FIhat1}), (\ref{gz.hat}), and the fact that $|\widehat{m}_m({\bf X}_i; \varphi)-Y_i|\leq 3L$, we can write
\begin{eqnarray}
S_n(2) &\leq & \big|\mathcal{F}_{\varepsilon}\big|\sup_{\varphi \in \mathcal{F}_{\varepsilon}}  P\Bigg\{ \ell^{-1} \sum_{i\in\,\boldsymbol{{\cal I}}_\ell} \left|\frac{1}{\widehat{\pi}_{\varphi}({\bf Z}_i, Y_i)}-\frac{1}{\pi_{\varphi}({\bf Z}_i, Y_i)}\right|  > \frac{\beta}{18L^2}\Bigg\} \nonumber\\
&\leq& 
\big|\mathcal{F}_{\varepsilon}\big|\sup_{\varphi \in \mathcal{F}_{\varepsilon}} \sum_{i\in\boldsymbol{{\cal I}}_\ell}E\Bigg[P\Bigg\{\Bigg|
\frac{1-\widetilde{\eta}_{m}({\bf Z}_i)}{\widetilde{\psi}_{m}({\bf Z}_i;\varphi)}
- \frac{E[1-\Delta_i|{\bf Z}_i]}{E[\Delta_i\,\varphi(Y_i)|{\bf Z}_i]}
\Bigg|\,\varphi(Y_i) >  \frac{\beta}{18L^2}\,\Bigg|{\bf Z}_i, Y_i\Bigg\}\Bigg],~~~~ \label{EQ7}
\end{eqnarray}
where the last line follows upon replacing the term $\exp\{g({\bf x})\}$ in the definition of $\pi_{\varphi}({\bf x},y)$ in (\ref{NonIg4}) by the right side of (\ref{expgz2}). Now, to bound (\ref{EQ7}), 
we note that  
\begin{eqnarray*}
	&& \left|\frac{1-\widetilde{\eta}_{m}({\bf Z}_i)}{ \widetilde{\psi}_{m}({\bf Z}_i; \varphi)} -
	\frac{1-E\big[\Delta_i\big|{\bf Z}_i\big]}{E\big[\Delta_i \varphi(Y_i)\big|{\bf Z}_i\big] }\right|\\
	&&~~~~~~~=~  \left|-\,\frac{1-\widetilde{\eta}_{m}({\bf Z}_i)}{ \widetilde{\psi}_{m}({\bf Z}_i; \varphi)} 
	\cdot \frac{\widetilde{\psi}_{m}({\bf Z}_i; \varphi) - E\big[\Delta_i\varphi(Y_i)\big|{\bf Z}_i\big]}{E\big[\Delta_i\varphi(Y_i)\big|{\bf Z}_i\big]}
	+ \frac{E\big[\Delta_i\big|{\bf Z}_i\big] -\widetilde{\eta}_{m}({\bf Z}_i) }{E\big[\Delta_i\varphi(Y_i)\big|{\bf Z}_i\big]}
	\right|\\
	&&~~~~~~~~~~~~~~~~\leq~
	\left|\frac{1-\widetilde{\eta}_{m}({\bf Z}_i)}{ \widetilde{\psi}_{m}({\bf Z}_i; \varphi)} \right|
	\cdot  \left|\frac{ \widetilde{\psi}_{m}({\bf Z}_i; \varphi) - E\big[\Delta_i\varphi(Y_i)\big|{\bf Z}_i\big]}{E\big[\Delta_i\varphi(Y_i)\big|{\bf Z}_i\big]}\right| 
	+ \left| \frac{E\big[\Delta_i\big|{\bf Z}_i\big] -\widetilde{\eta}_{m}({\bf Z}_i) }{E\big[\Delta_i\varphi(Y_i)\big|{\bf Z}_i\big]}\right|.
\end{eqnarray*}
Therefore, in view of (\ref{psieta}), the inner conditional probability in (\ref{EQ7}) becomes
\begin{eqnarray}
&& P\left\{\left|\frac{1-\widetilde{\eta}_{m}({\bf Z}_i)}{ \widetilde{\psi}_{m}({\bf Z}_i; \varphi)} -
\frac{1-E\big[\Delta_i\big|{\bf Z}_i\big]}{E\big[\Delta_i \varphi(Y_i)\big|{\bf Z}_i\big] }\right|\, \varphi(Y_i)\, >\,\frac{\beta}{18L^2}\,\Bigg| {\bf Z}_i, Y_i\right\} \nonumber\\
&& ~~~~~~\leq~  P\left\{\left| \frac{1-\widetilde{\eta}_{m}({\bf Z}_i)}{\widetilde{\psi}_{m}({\bf Z}_i; \varphi)}\right|\cdot  \left|\widetilde{\psi}_{m}({\bf Z}_i; \varphi) - E\big[\Delta_i\varphi(Y_i)\big|{\bf Z}_i\big]\right| > \frac{\varrho_0\beta}{36B L^2} \bigg| {\bf Z}_i, Y_i\right\} ~~~~~~\nonumber\\
&& ~~~~~~~~~~+\,  P\left\{\Big|\widetilde{\eta}_{m}({\bf Z}_i) - E\big[\Delta_i\big|{\bf Z}_i\big]
 \Big| > \frac{\varrho_0\beta}{36B L^2} \bigg| {\bf Z}_i, Y_i\right\}\nonumber \\[3pt]
&& ~~~~~~:=~ P_{n,1}(i) +P_{n,2}(i),     \label{EEQQ1}
\end{eqnarray}
where we have used the facts that $\varphi(y)\in (0, B]$, $B>0$, and $E\big[\Delta \varphi(Y)\big|{\bf Z}={\bf z}\big] \geq \varrho_0$ (in view of by Assumption D). But, using standard arguments, it is not difficult to show that, under assumptions (B)--(E) and  $m$ large enough, one has 
\begin{equation}\label{EEQQ2}
P_{n,2}(i) ~\leq~ C_{12} \,e^{-C_{13}\,m h^d \beta^2}
\end{equation}
where $c_{12}$ and $c_{13}$ are positive constants not depending om $m$, $\ell$, or $\beta$. Next, to deal with the term $P_{n,1}(i)$, put
\[
\mathcal{B}_m({\bf Z}_i) = \big\{ \widetilde{\psi}_{m}({\bf Z}_i; \varphi) \geq \varrho_0/2\big\},
\]
where $\varrho_0$ is as in assumption D, and note that
\begin{eqnarray*}
P_{n,1}(i) &\leq&  P\Bigg\{\Bigg[\left| \frac{1-\widetilde{\eta}_{m}({\bf Z}_i)}{\widetilde{\psi}_{m}({\bf Z}_i; \varphi)}\right|\cdot  \left|\widetilde{\psi}_{m}({\bf Z}_i; \varphi) - E\big[\Delta_i\varphi(Y_i)\big|{\bf Z}_i\big]\right| > \frac{\varrho_0\beta}{36B L^2} \Bigg] \cap \mathcal{B}_m({\bf Z}_i)\bigg| {\bf Z}_i, Y_i\Bigg\}\\[3pt]
&&~ +~ P\Big\{\mathcal{B}^c_m({\bf Z}_i) \big| {\bf Z}_i, Y_i \Big\} \\[3pt]
&:=& P'_{n,1}(i) + P''_{n,1}(i).
\end{eqnarray*}
However, straightforward but tedious arguments show that
\[
 P'_{n,1}(i) ~\leq~  P\bigg\{ \left|\widetilde{\psi}_{m}({\bf Z}_i; \varphi) - E\big[\Delta_i\varphi(Y_i)\big|{\bf Z}_i\big]\right| > \frac{\varrho^2_0\beta}{72 B L^2}\bigg| {\bf Z}_i, Y_i\bigg\} ~\leq~ C_{14}\,e^{-C_{15}\,m h^d \beta^2},
\]
for $n$ (and thus $m$) large enough, where $C_{14}$ and $C_{15}$ are positive constants not depending on $m$, $\ell$, or $\beta$. 
As for the term $P''_{n,1}(i)$, we have
$
P''_{n,1}(i) =  P\big\{\widetilde{\psi}_{m}({\bf Z}_i; \varphi) - E\big[\Delta_i\varphi(Y_i)\big|{\bf Z}_i\big] < \varrho_0/2 - E\big[\Delta_i\varphi(Y_i)\big|{\bf Z}_i\big]\big| {\bf Z}_i, Y_i\big\} \leq   
P\big\{\big|\widetilde{\psi}_{m}({\bf Z}_i; \varphi) - E\big[\Delta_i\varphi(Y_i)\big|{\bf Z}_i\big] \big| > 
\varrho_0/2\,\big| {\bf Z}_i, Y_i\big\}\,\leq C_{16}\,
\exp\{-C_{17}\,m h^d\},
$
where we have used the fact that $\psi_2$ is bounded by assumption (D); here $C_{16}$ and  $C_{17}$ are positive constant not depending on $m$ or $\ell$. Putting these bounds together, we find
\begin{equation}\label{EEQQ3}
P_{n,1}(i) \, \leq P'_{n,1}(i) + P''_{n,1}(i) ~\leq\, C_{14}\,e^{-C_{15}\,m h^d \beta^2} + C_{16}\, e^{-C_{17}\,m h^d}.
\end{equation}
Therefore, in view of (\ref{EQ5}) -- (\ref{EEQQ3}), for every $\beta>0$ and $n$ large enough, we have
\begin{eqnarray}
P\left\{ \sup_{\varphi \in \mathcal{F}_{\varepsilon}} \,\left| \widehat{L}_{m,\ell}(\varphi) - E\left[\Big|\widehat{m}_m({\bf X}; \varphi) - Y\Big|^2\Big|\mathbb{D}_m\right] \right|  > \beta \right\} &\leq &  \ell\,\big|\mathcal{F}_\varepsilon\big| \bigg[C_{21}\,e^{-C_{22} mh^d} +C_{23}\, e^{-C_{24}\,m h^d \beta^2}\bigg]\nonumber\\&&~~~~+\,2\big|\mathcal{F}_\varepsilon\big|\, 
e^{ -\pi^2_{\mbox{\tiny min}} \ell \beta^2/(162 L^4)},  \label{ADD1}
\end{eqnarray}
Next, we deal with the second term on the right side of (\ref{EQNN}). To this end,  first note that by (\ref{STAR1}) and the fact that $E\big|m({\bf X}; \varphi) - Y\big|^2 = E\big[\Delta\big|m({\bf X}; \varphi) - Y\big|^2 \big/\pi_{\varphi}({\bf Z}, Y)\big]$, one has
\begin{eqnarray}
\left| \widehat{L}_{m,\ell}(\varphi) - E\Big|m({\bf X}; \varphi) - Y\Big|^2 \right| &\leq& 
\left|\frac{1}{\ell}\, \sum_{i\in\,\boldsymbol{{\cal I}}_\ell} \frac{\Delta_i \,\big|\widehat{m}_{m}({\bf X}_i;\varphi)-Y_i\big|^2}{ \pi_{\varphi}({\bf Z}_i, Y_i)}
- \frac{1}{\ell}\, \sum_{i\in\,\boldsymbol{{\cal I}}_\ell} \frac{\Delta_i \,\big|m({\bf X}_i;\varphi)-Y_i\big|^2}{ \pi_{\varphi}({\bf Z}_i, Y_i)}\right| \nonumber\\
&& +\,\left| \frac{1}{\ell}\, \sum_{i\in\,\boldsymbol{{\cal I}}_\ell} \frac{\Delta_i \,\big|m({\bf X}_i;\varphi)-Y_i\big|^2}{ \pi_{\varphi}({\bf Z}_i, Y_i)}
- E\left[\frac{\Delta\big|m({\bf X}; \varphi) - Y\big|^2}{\pi_{\varphi}({\bf Z}, Y)}\right]\right| \nonumber\\
&& ~+\, \left|\frac{1}{\ell} \sum_{i\in\,\boldsymbol{{\cal I}}_\ell} \Delta_i \,
\big|\widehat{m}_{m}({\bf X}_i;\varphi)-Y_i\big|^2 \left[\frac{1}{\pi_{\varphi}({\bf Z}_i, Y_i)} -\frac{1}{\widehat{\pi}_{\varphi}({\bf Z}_i, Y_i)}\right]\right| \nonumber \\[3pt]
&:=& |U_{n,1}(\varphi)| +  |U_{n,2}(\varphi)| +  |U_{n,3}(\varphi)|. \label{EQS1}
\end{eqnarray}
Therefore, for every $\beta>0$, 
\begin{eqnarray*}
P\left\{\sup_{\varphi \in \mathcal{F}_{\varepsilon}}\left| \widehat{L}_{m,\ell}(\varphi) - E\Big|m({\bf X}; \varphi) - Y\Big|^2 \right| > \beta\right\} &\leq&
\sum_{k=1}^3 P\left\{\sup_{\varphi \in \mathcal{F}_{\varepsilon}}\big| U_{n,k}(\varphi)\big| > \frac{\beta}{3}\right\}.
\end{eqnarray*}
But using Assumption (A) and the simple fact that $a^2-b^2\leq |a-b||a+b|$, one can write
\begin{eqnarray*}
&& P\bigg\{\sup_{\varphi \in \mathcal{F}_{\varepsilon}}\big| U_{n,1}(\varphi)\big| > \beta/3\bigg\}\\
&&\leq\,
P\bigg\{\sup_{\varphi \in \mathcal{F}_{\varepsilon}}\frac{1}{\ell} \sum_{i\in\,\boldsymbol{{\cal I}}_\ell} \Big[\big|\widehat{m}_{m}({\bf X}_i;\varphi)-
m({\bf X}_i; \varphi)\big|\cdot \big|\widehat{m}_{m}({\bf X}_i;\varphi)+m({\bf X}_i; \varphi) -2Y_i\big|\Big] > \frac{\beta \pi_{\mbox{\tiny min}}}{3}\bigg\}\\
&&\leq\,
\big|\mathcal{F}_{\varepsilon}\big| \sup_{\varphi \in \mathcal{F}_{\varepsilon}} \sum_{i\in\,\boldsymbol{{\cal I}}_\ell} 
P\left\{\big|\widehat{m}_{m}({\bf X}_i;\varphi) - m({\bf X}_i; \varphi)\big| >\frac{\beta \pi_{\mbox{\tiny min}}}{15 L}\right\},
\end{eqnarray*}
where we have used the fact that $\big|\widehat{m}_{m}({\bf X}_i;\varphi)+m({\bf X}_i; \varphi) -2Y_i\big|\leq 5L$. Now,  using standard arguments, it is not hard to show that under assumptions (B)--(E) and  $m$ large enough, one has
\begin{equation}\label{EQS2}
P\Big\{\sup_{\varphi \in \mathcal{F}_{\varepsilon}}| U_{n,1}(\varphi)| > \beta/3\Big\} \,\leq\, \ell\,\big|\mathcal{F}_{\varepsilon}\big|C_{25}\, \exp\{-C_{26}\,m h^d \beta^2\},
\end{equation}
for positive constants $C_{25}$ and $C_{26}$ not depending on $m$, $\ell$, or $\beta$. Next, since the iid random variabels $\Delta_i |m({\bf X}_i; \varphi) - Y_i|^2/\pi_{\varphi}({\bf Z}_i, Y_i)$, $i\in\,\boldsymbol{{\cal I}}_\ell$, take values in $(0,\,4L^2/\pi_{\mbox{\tiny min}})$, an application of  Hoeffding's inequality yields 
\begin{equation}\label{EQS3}
P\Big\{\sup_{\varphi \in \mathcal{F}_{\varepsilon}}| U_{n,2}(\varphi)| > \beta/3\Big\}\,\leq\, \big|\mathcal{F}_{\varepsilon}\big| \sup_{\varphi \in \mathcal{F}_{\varepsilon}} P\{| U_{n,2}(\varphi)| > \beta/3\} \,\leq\, 2\,\big|\mathcal{F}_{\varepsilon}\big| \exp\big\{-\ell \pi^2_{\mbox{\tiny min}}\beta^2/(72 L^4)\big\}.
\end{equation}
Finally, the same arguments that were used to deal with the term $S_n(2)$ in  (\ref{EQ5}) can be employed to show that
\begin{equation}\label{EQS4}
P\Big\{\sup_{\varphi \in \mathcal{F}_{\varepsilon}}| U_{n,3}(\varphi)| > \beta/3\Big\}\,\leq\, \ell \, \big|\mathcal{F}_{\varepsilon}\big| \left[ C_{27}\,\exp\big\{-C_{28}\,m h^d\big\} + C_{29}\, \exp\big\{-C_{30}\,m h^d \beta^2\big\}\right],
\end{equation}
for $n$ large enough and positive constants $C_{27}$ -- $C_{30}$  that do not depend on $m$, $\ell$, or $\beta$. Putting together (\ref{EQS1}), (\ref{EQS2}), (\ref{EQS3}), and (\ref{EQS4}), one finds, for every$\beta>0$,
\begin{eqnarray}
P\left\{\sup_{\varphi \in \mathcal{F}_{\varepsilon}}\left| \widehat{L}_{m,\ell}(\varphi) - E\Big|m({\bf X}; \varphi) - Y\Big|^2 \right| >\beta\right\}&\leq&
C_{31}\,\ell\,\big|\mathcal{F}_{\varepsilon}\big| \left[e^{-C_{32}\,m h^d} + e^{-C_{33}\,m h^d \beta^2}\right]
\nonumber\\  &&~~~~~~~~
+\,2\,\big|\mathcal{F}_{\varepsilon}\big| \,e^{-C_{34}\,\ell\, \beta^2},   \label{EQS5}
\end{eqnarray}
for $n$ large enough, where $C_{31}$ -- $C_{34}$ are positive constants not depending on $m$ or $\ell$. Now to complete the proof of the theorem, let $0<\varepsilon_n \downarrow 0$ be as in the statement of the theorem and let $\varphi_{\varepsilon_n}$ be as in (\ref{FIEP}). Then, in view of (\ref{repr2}),  (\ref{EQ4B}), (\ref{EEQ2}), and the $C_p$-inequality (with $p=2$), one has
\begin{eqnarray}
\int\Big|\widehat{m}({\bf x}; \widehat{\varphi}_n)- m({\bf x}) \Big|^2  \mu(d{\bf x}) 
&\leq& 
2\int\Big|\widehat{m}({\bf x}; \widehat{\varphi}_n)- m({\bf x}; \varphi_{\varepsilon_n}) \Big|^2  \mu(d{\bf x})\,+ \,4LC \,\varepsilon_n\,, \label{ABC2}
\end{eqnarray}
where $C>0$ is the constant in Lemma \ref{LEM-2}. Now observe that  (\ref{ABC2}) in conjunction with Lemma \ref{LEM-1} implies that, for every $t>0$,
\begin{eqnarray*}
&& P\left\{\int\Big|\widehat{m}({\bf x}; \widehat{\varphi}_n)- m({\bf x}) \Big|^2  \mu(d{\bf x})  >   t\right\} 
\leq
P\left\{\int\Big|\widehat{m}({\bf x}; \widehat{\varphi}_n)- m({\bf x}; \varphi_{\varepsilon_n}) \Big|^2  \mu(d{\bf x}) > \frac{t}{2}-2LC\, \varepsilon_n\right\}\\
&&~ \leq~ 
P\left\{\sup_{\varphi \in \mathcal{F}_{\varepsilon_n}} \,\Bigg|E\left[\Big|\widehat{m}_m({\bf X}; \varphi) - Y\Big|^2\Big|\mathbb{D}_m\right] - 
\widehat{L}_{m,\ell}(\varphi) \Bigg| \,>\,  \frac{t/2 - 2LC\, \varepsilon_n -C_1\sqrt{\varepsilon_n}}{2}\right\}\\
&&~~~~~~+\, P\left\{\sup_{\varphi \in \mathcal{F}_{\varepsilon_n}} \,\Bigg|\widehat{L}_{m,\ell}(\varphi) - E\Big|m({\bf X}; \varphi) - Y\Big|^2 \Bigg| \,>\,  \frac{t/2 - 2LC\, \varepsilon_n -C_1\sqrt{\varepsilon_n}}{2}\right\}.
\end{eqnarray*}
Now, since $\varepsilon_n\downarrow 0$, as $n\to\infty$, we can choose $n$  large enough so that $t/2-2LC\,\varepsilon_n - C_1 \sqrt{\varepsilon_n}\,>t/4$. Therefore, in view of (\ref{ADD1}) and  (\ref{EQS5}), for every $t>0$ and for $n$ large enough, one finds
\begin{equation*}
P\left\{\int\Big|\widehat{m}({\bf x; \widehat{\varphi}_n})- m({\bf x}) \Big|^2  \mu(d{\bf x})  >   t\right\} 
\,\leq\,
4\,\big|\mathcal{F}_{\varepsilon_n}\big|\, e^{-c_{35}\ell t^2} + c_{36}\, \ell \big|\mathcal{F}_{\varepsilon_n}\big|\, e^{-c_{37}m h^d} + c_{38}\, \ell \big|\mathcal{F}_{\varepsilon_n}\big|\, e^{-c_{39}m h^d t^2},
\end{equation*}
which completes the proof of  Theorem \ref{THM-B}.

\hfill $\Box$

\vspace{3.5mm}\noindent
PROOF OF COROLLARY \ref{COR-A}

\vspace{1mm}\noindent
Corollary \ref{COR-A} follows from an application of the Borel-Cantelli lemma in conjunction with (\ref{cond}), the bound in Theorem \ref{THM-B}, and Remark \ref{REM-thm2}.

\hfill $\Box$

\vspace{3.5mm}\noindent
PROOF OF THEOREM \ref{THM-BB}

\vspace{1mm}\noindent
We first note that, by Remark \ref{REM-thm2}, it is sufficient to prove the theorem for the case of $p=2$.  The proof is along standard arguments and goes as follows. Observe that
\begin{eqnarray}
E\big|\widehat{m}({\bf X}; \widehat{\varphi}_n)- m({\bf X}) \big|^2
&=& E\left[ \int_{\mathbb{R}^d}\Big|\widehat{m}({\bf x}; \widehat{\varphi}_n)- m({\bf x}) \Big|^2\,\mu(d{\bf x}  \right]\nonumber\\
&=& \int_0^{\infty}P\left\{\int_{\mathbb{R}^d}\Big|\widehat{m}({\bf x}; \widehat{\varphi}_n)- m({\bf x}) \Big|^2  \mu(d{\bf x}) >  t \right\} dt\nonumber\\
&=& \int_0^{9L^2}P\left\{\int_{\mathbb{R}^d}\Big|\widehat{m}({\bf x}; \widehat{\varphi}_n)- m({\bf x}) \Big|^2  \mu(d{\bf x}) >  t \right\} dt\,, \label{B9L2}
\end{eqnarray}
where the last line follows from the fact that, by the definition of the estimator $\widehat{m}({\bf x}; \widehat{\varphi}_{n})$ in (\ref{mhat5}), 
\begin{eqnarray*}
\Big|\widehat{m}({\bf x}; \widehat{\varphi}_n)- m({\bf x}) \Big|^2
&\leq& \big(\big|\widehat{m}({\bf x}; \widehat{\varphi}_n)\big|+\big|m({\bf x})\big|\big)^2 \\
&\leq& \left(\big|\widehat{\eta}_{m,1}({\bf x})\big|+\left|\frac{\widehat{\psi}_{m,1}({\bf x}; \widehat{\varphi}_n)}{\widehat{\psi}_{m,2}({\bf x}; \widehat{\varphi}_n)}\right|\cdot\big|1-\widehat{\eta}_{m,2}({\bf x})\big| +L	\right)^2 
\,\leq~ (L+L\cdot 1 +L)^2.
\end{eqnarray*}
Therefore, by Theorem \ref{THM-B}, for $n$ large enough, we have
\begin{eqnarray}
&& \big(\mbox{right side of (\ref{B9L2})}\big)\nonumber\\
&&~\leq~ \int_0^u dt \,+\, \big(c_4\vee c_6\big) \big|\mathcal{F}_{\varepsilon_n}\big|\cdot\Bigg[
\int_u^{9L^2} e^{-c_5\ell t^2}\, dt \,+\,\ell \int_u^{9L^2} e^{-c_8\,mh^d t^2}\, dt \, +\,\ell\, e^{-c_7\,mh^d} \int_u^{9L^2} dt\,
\Bigg], \nonumber\\[4pt]
&& ~~~~~~~~~\mbox{(where $c_4$--$c_8$ are as in Theorem \ref{THM-B})} \nonumber\\
&&~\leq~ u\,+\, 2(c_4\vee c_6) \big|\mathcal{F}_{\varepsilon_n}\big|\, \ell
\int_u^{9L^2} e^{-(c_5\wedge c_8)(\ell\wedge mh^d)\, t^2}\, dt \, +
(c_4 \vee c_6)(9 L^2) \big|\mathcal{F}_{\varepsilon_n}\big|\, \ell\,e^{-c_7\,mh^d} \nonumber\\[4pt]
&&~\leq~ u\,+\, \frac{2(c_4\vee c_6) \big|\mathcal{F}_{\varepsilon_n} \big|\, \ell}{\sqrt{(c_5\wedge c_8)(\ell \wedge mh^d)}}\cdot
\int_{u\sqrt{(c_5\wedge c_8)(\ell \wedge mh^d)}}^{\infty}\, e^{-v^2/2}\, dv \, +\, (c_4 \vee c_6)(9 L^2) \big|\mathcal{F}_{\varepsilon_n}\big|\,\ell\, e^{-c_7\,mh^d} \nonumber\\[4pt]
&& ~~~~~~~~~\mbox{(which follows from the change of variable $v=\sqrt{(c_5\wedge c_8)(\ell \wedge mh^d)} \,\, t$\,} \nonumber\\[4pt]
&& ~\leq~ u\,+\, \frac{2(c_4\vee c_6) \big|\mathcal{F}_{\varepsilon_n} \big|\, \ell}{\sqrt{(c_5\wedge c_8)(\ell \wedge mh^d)}}\cdot \frac{e^{-(c_5\wedge c_8)(\ell \wedge mh^d)\,u^2/2}}{\sqrt{(c_5\wedge c_8)(\ell \wedge mh^d)} \,\, u}
\, +\, (c_4 \vee c_6)(9 L^2) \big|\mathcal{F}_{\varepsilon_n}\big|\,\ell\, e^{-c_7\,mh^d}, \label{mills}
\end{eqnarray}
where the last line follows from the upper bound on Mill's ratio; see, for example, Mitrinovic (1970; p. 177). Now, put 
\[
c=2(c_4 \vee c_6) \big|\mathcal{F}_{\varepsilon_n}\big|\,\ell~~~~\mbox{and}~~~
N=(c_5\wedge c_8)(\ell \wedge mh^d)/4
\]
and observe that the right side of (\ref{mills}) can be written as
\begin{equation}
 u + \frac{c}{4Nu}\,e^{-2Nu^2}+(c_4 \vee c_6)(9 L^2) \big|\mathcal{F}_{\varepsilon_n}\big|\,\ell\, e^{-c_7\,mh^d}. \label{CEQ9}
\end{equation}
But the term  $u+\frac{c}{4Nu}\,e^{-2Nu^2}$ in (\ref{CEQ9}) is approximately minimized by taking $u=\sqrt{\log (c)/ (2N)}$, and the corresponding minimum value of (\ref{CEQ9}) is
\begin{eqnarray*}
&& \sqrt{\frac{\log(c)}{2N}} + \sqrt{\frac{1}{8N\log(c)}}+\, (c_4 \vee c_6)(9 L^2) \big|\mathcal{F}_{\varepsilon_n}\big|\,\ell\, e^{-c_7\,mh^d}\\
&&~~=~ \sqrt{\frac{a_1+\log \ell+ \log|\mathcal{F}_{\varepsilon_n}|}{a_2\,(\ell \wedge mh^d)}}
\,+ \sqrt{\frac{1}{a_3\,(\ell \wedge mh^d)\big[ a_1+\log \ell+\log|\mathcal{F}_{\varepsilon_n}|\big]}}
\,\,+ a_4\,\big|\mathcal{F}_{\varepsilon_n}\big|\,\ell\, e^{-c_7\,mh^d},~~
\end{eqnarray*}
where $a_1$--\,$a_4$ are positive constants not depending on $m$, $\ell$, and $n$.

\hfill $\Box$

\vspace{4mm}\noindent
PROOF OF THEOREM \ref{THM-BBC}

\vspace{1mm}\noindent
Let $\widehat{m}^{\mbox{\tiny HT}}_{m}({\bf x};\widetilde{\pi}_\varphi)$,  $m({\bf x}, \pi_{\varphi^*})$ , and $\varphi_{\varepsilon}$ be as in    (\ref{CEQ4}), 
(\ref{FS15}), and   (\ref{FIEP2}) respectively.
Also, define 
\begin{equation}\label{UA2}
\widetilde{L}_{m,\ell}(\widetilde{\pi}_{\varphi}) = 
\ell^{-1} \sum_{i\in\,\boldsymbol{{\cal I}}_\ell} \frac{\Delta_i}{ \widetilde{\pi}_{\varphi}({\bf Z}_i, Y_i)}\,\Big| 
\widehat{m}^{\mbox{\tiny HT}}_{m}({\bf  X}_i; \widetilde{\pi}_\varphi) - Y_i\Big|^2 ,
\end{equation}
where $\widetilde{\pi}_{\varphi}({\bf x}, y)$  is given by (\ref{CEQ1}), and put
\[
\widetilde{\varphi}_{\varepsilon} = \argmin_{\varphi \in \mathcal{F}_{\varepsilon}}\widetilde{L}_{m,\ell}(\widetilde{\pi}_{\varphi}).
\]
Then, using the arguments that led to (\ref{EQ4}) and (\ref{EQ4Q}), yield
\begin{eqnarray}
&& \int\Big|\widehat{m}^{\mbox{\tiny HT}}_m({\bf x}; \widetilde{\pi}_{\widetilde{\varphi}_{\varepsilon}}) - m({\bf x};\pi_{\varphi_{\varepsilon}})\Big|^2  \mu(d{\bf x})\nonumber\\
&& \leq \sup_{\varphi \in \mathcal{F}_{\varepsilon}} \,\left|E\left[\Big|\widehat{m}^{\mbox{\tiny HT}}_m({\bf X}; \widetilde{\pi}_{\varphi}) - Y\Big|^2\Big|\mathbb{D}_m\right] - \widetilde{L}_{m,\ell}(\widetilde{\pi}_{\varphi}) \right|  +  \sup_{\varphi \in \mathcal{F}_{\varepsilon}} \,\left| \widetilde{L}_{m,\ell}(\widetilde{\pi}_{\varphi}) - 
E\Big|m({\bf X}; \pi_{\varphi}) - Y\Big|^2 \right|\nonumber\\
&&~~~~ +~ 2\,E\left[\big|\widehat{m}^{\mbox{\tiny HT}}_m({\bf X}; \widetilde{\pi}_{\widetilde{\varphi}_{\varepsilon}}) - m({\bf X};\pi_{\varphi_{\varepsilon}})\big|\cdot
\big| m({\bf X};\pi_{\varphi_{\varepsilon}}) -m({\bf X};\pi_{\varphi^{*}}) \big| \Big|\mathbb{D}_n  \right] \label{CSC}
\end{eqnarray}
where, as before,  $\varphi^*$ is the true $\varphi$. But by Cauchy-Schwarz inequality, the last line on the right side of (\ref{CSC}) is bounded by 
\begin{eqnarray}
&& 2\,\sqrt{\int\Big|\widehat{m}^{\mbox{\tiny HT}}_m({\bf x}; \widetilde{\pi}_{\widetilde{\varphi}_{\varepsilon}}) - m({\bf x};\pi_{\varphi_{\varepsilon}})\Big|^2  \mu(d{\bf x})}\cdot \sqrt{E\big| m({\bf X};\pi_{\varphi_{\varepsilon}}) -m({\bf X};\pi_{\varphi^{*}}) \big|^2} \nonumber\\
&&~\leq ~ C_3\, \sqrt{\int\Big|\widehat{m}^{\mbox{\tiny HT}}_m({\bf x}; \widetilde{\pi}_{\widetilde{\varphi}_{\varepsilon}}) - m({\bf x};\pi_{\varphi_{\varepsilon}})\Big|^2  \mu(d{\bf x})}\, \cdot\sqrt{\varepsilon\,}\,,\label{CSC2}
\end{eqnarray}
where (\ref{CSC2}) follows from arguments similar to those used to arrive at (\ref{EQ4B}) and (\ref{EEQ2}); here $C_3$ is a positive constant not depending on $n$ or $\varepsilon$. Therefore, in view of  (\ref{CSC}) and (\ref{CSC2}), for any $t>0$
\begin{eqnarray}
P\left\{
\int\Big|\widehat{m}^{\mbox{\tiny HT}}_m({\bf x}; \widetilde{\pi}_{\widetilde{\varphi}_{\varepsilon}}) - m({\bf x};\pi_{\varphi_{\varepsilon}})\Big|^2  \mu(d{\bf x}) > t
\right\} - P\left\{\int\Big|\widehat{m}^{\mbox{\tiny HT}}_m({\bf x}; \widetilde{\pi}_{\widetilde{\varphi}_{\varepsilon}}) - 
m({\bf x};\pi_{\varphi_{\varepsilon}})\Big|^2  \mu(d{\bf x}) > \frac{t^2}{c_4 \varepsilon} \right\}  \nonumber \\[3pt]
\leq\, P\left\{ \sup_{\varphi \in \mathcal{F}_{\varepsilon}} \,\left|E\left[\Big|\widehat{m}^{\mbox{\tiny HT}}_m({\bf X}; \widetilde{\pi}_{\varphi}) - Y\Big|^2\Big|\mathbb{D}_m\right] - \widetilde{L}_{m,\ell}(\widetilde{\pi}_{\varphi}) \right|>\frac{t}{3}\right\}~~~~~\nonumber \\[3pt]
+\, 
P\left\{
\sup_{\varphi \in \mathcal{F}_{\varepsilon}} \,\left| \widetilde{L}_{m,\ell}(\widetilde{\pi}_{\varphi}) - 
E\Big|m({\bf X}; \pi_{\varphi}) - Y\Big|^2 \right|>\frac{t}{3}\right\},~~~~~~~~~~~~~ \label{Bound99}
\end{eqnarray}
where $c_4=(3C_3)^2$ with $C_3$ as in (\ref{CSC2}). But observe that for every constant $\beta>0$
\begin{eqnarray}
&& P\left\{ \sup_{\varphi \in \mathcal{F}_{\varepsilon}} \,\left| \widetilde{L}_{m,\ell}(\widetilde{\pi}_\varphi) - E\left[\Big|\widehat{m}^{\mbox{\tiny HT}}_m({\bf X}; \widetilde{\pi}_{\varphi}) - Y\Big|^2\Big|\mathbb{D}_m\right] \right|  > \beta \right\} \nonumber\\
&&~\leq\,   P\left\{ \sup_{\varphi \in \mathcal{F}_{\varepsilon}} \,\left|\ell^{-1} \sum_{i\in\,\boldsymbol{{\cal I}}_\ell} \frac{\Delta_i \,
\big|\widehat{m}^{\mbox{\tiny HT}}_{m}({\bf X}_i;\widetilde{\pi}_{\varphi})-Y_i\big|^2}{\pi_{\varphi}({\bf Z}_i, Y_i)}  - E\left[\Big|\widehat{m}^{\mbox{\tiny HT}}_m({\bf X}; \widetilde{\pi}_{\varphi}) - Y\Big|^2\Big|\mathbb{D}_m\right] \right| > \frac{\beta}{2}\right\} \nonumber\\
&&~~~~~~~+\, P\left\{ \sup_{\varphi \in \mathcal{F}_{\varepsilon}} \,\left|\ell^{-1} \sum_{i\in\,\boldsymbol{{\cal I}}_\ell} \Delta_i \,
\big|\widehat{m}^{\mbox{\tiny HT}}_{m}({\bf X}_i; \widetilde{\pi}_{\varphi})-Y_i\big|^2 \left[\Big(\pi_{\varphi}({\bf Z}_i, Y_i)\Big)^{-1}-\Big(\widetilde{\pi}_{\varphi}({\bf Z}_i, Y_i)\Big)^{-1}\right]\right|>\frac{\beta}{2}\right\} \nonumber\\[4pt]
&& ~:=\, T_n(1)+T_n(2).  \label{Tn12}
\end{eqnarray}
On the other hand, for every $i\in \boldsymbol{{\cal I}}_\ell\,,$  we find that 
$
E\big[  \Delta_i |\widehat{m}^{\mbox{\tiny HT}}_m({\bf X}_i; \widetilde{\pi}_{\varphi}) -Y_i |^2 / \pi_{\varphi}({\bf Z}_i, Y_i)\,\big| \mathbb{D}_m\big] =
E\big[ E\big\{ \Delta_i |\widehat{m}^{\mbox{\tiny HT}}_m({\bf X}_i; \widetilde{\pi}_{\varphi}) -Y_i |^2 / \pi_{\varphi}({\bf Z}_i, Y_i)\,\big| \mathbb{D}_m, {\bf X}_i, Y_i\big\}\big| \mathbb{D}_m\big] = E\big[ |\widehat{m}^{\mbox{\tiny HT}}_m({\bf X}_i; \widetilde{\pi}_{\varphi}) -Y_i |^2 \,\big| \mathbb{D}_m\big],$ where the last expression follows from the definition of $\pi_{\varphi}$ in (\ref{NonIg4}). Moreover, by the definition of $\widetilde{\pi}_{\varphi}({\bf Z}, Y)$, as given by (\ref{CEQ1}), one finds
\begin{eqnarray}
|\widehat{m}^{\mbox{\tiny HT}}_m({\bf X}_i; \widetilde{\pi}_{\varphi})| &\leq& \max_{k\in \boldsymbol{{\cal I}}_m}
\big|\Delta_k Y_k / \widetilde{\pi}_{\varphi}({\bf Z}_k, Y_k)\big| 
~\leq~ 
L\cdot\left( 1+   \max_{k\in \boldsymbol{{\cal I}}_m}\left|\frac{1}{\widetilde{\psi}_m({\bf Z}_k; \varphi)}\right|\cdot B  \right),   \label{OABC}
\end{eqnarray}
where the function $\widetilde{\psi}_m({\bf Z}_k;\varphi)$ is as given in (\ref{ETPS}). Consequently, conditional on $\mathbb{D}_m$, the terms $\big[\Delta_i |\widehat{m}^{\mbox{\tiny HT}}_m({\bf X}_i; \widetilde{\pi}_{\varphi}) -Y_i |^2\big] \big/ \pi_{\varphi}({\bf Z}_i, Y_i)$, $i\in \boldsymbol{{\cal I}}_{\ell}$,\, are independent  nonnegative random variables bounded by $2L^2\big\{ 4+B^2 \max^2_{k\in \boldsymbol{{\cal I}}_m}\big|1/\widetilde{\psi}_m({\bf Z}_k; \varphi)\big|   \big\}/\pi_{\mbox{\tiny min}}.$ Therefore, the term $T_n(1)$ in (\ref{Tn12}) can be handled as follows
\begin{eqnarray}
&& T_n(1)    \nonumber\\
&&\leq  \big|\mathcal{F}_{\varepsilon}\big|
\sup_{\varphi \in \mathcal{F}_{\varepsilon}}  
E \Bigg[ P\Bigg\{ \Bigg|\ell^{-1} \sum_{i\in\boldsymbol{{\cal I}}_\ell} \frac{\Delta_i  \big|\widehat{m}^{\mbox{\tiny HT}}_{m}({\bf X}_i;\widetilde{\pi}_{\varphi})-Y_i\big|^2}{\pi_{\varphi}({\bf Z}_i, Y_i)}  -  
E\left[\Big|\widehat{m}^{\mbox{\tiny HT}}_m({\bf X}; \widetilde{\pi}_{\varphi}) - Y\Big|^2\Big|\mathbb{D}_m\right] \Bigg| > \frac{\beta}{2} \bigg| \mathbb{D}_m\Bigg\}\Bigg]  \nonumber\\
& &\leq 2\,\big|\mathcal{F}_{\varepsilon}\big|\,\sup_{\varphi \in \mathcal{F}_{\varepsilon}}  E\left[\exp\left\{  \frac{-2\ell^2\,(\beta/2)^2}{2L^2 \pi^{-1}_{\mbox{\tiny min}}  
\left\{ 4+B^2 \max^2_{k\in \boldsymbol{{\cal I}}_m}\big|1/\widetilde{\psi}_m({\bf Z}_k; \varphi)\big|   \right\}}\right\} \right],   \label{ESP1}
\end{eqnarray}
via Hoeffding's inequality. Now let $\varrho_0$ be the constant in Assumption (D) and observe that since the exponential function in (\ref{ESP1}) is always bounded by 1, the expectation on the right side of (\ref{ESP1}) is bounded by
\begin{eqnarray}
&& E\left[\exp\left\{  \frac{-2\ell^2\,(\beta/2)^2}{2L^2 \pi^{-1}_{\mbox{\tiny min}}  
\left\{ 4+\max^2_{k\in \boldsymbol{{\cal I}}_m}\big|B/\widetilde{\psi}_m({\bf Z}_k; \varphi)\big|   \right\}}\right\}\, \mathbb{I}\bigg\{\bigcap_{k\in\,\boldsymbol{{\cal I}}_m}\left[ \widetilde{\psi}_m({\bf Z}_k; \varphi) \geq\frac{\varrho_0}{2} \right]   \bigg\}\right] \nonumber\\[3pt]
&& ~~~~~~ +\,  E\Bigg[ \mathbb{I}\bigg\{\bigcup_{k\in\,\boldsymbol{{\cal I}}_m}\left[ \widetilde{\psi}_m({\bf Z}_k; \varphi) < \varrho_0/2 \right]   \bigg\}  \Bigg]\nonumber\\
&&~\leq~   \exp\left\{  \frac{-2\ell^2\,(\beta/2)^2}{2L^2 \pi^{-1}_{\mbox{\tiny min}}  
\big\{ 4+B^2 (2/\varrho_0)^2  \big\}}\right\}\cdot P\bigg\{\bigcap_{k\in\,\boldsymbol{{\cal I}}_m}\left[ \widetilde{\psi}_m({\bf Z}_k; \varphi) \geq\varrho_0/2 \right]   \bigg\}\ \nonumber\\[3pt]
&& ~~~~~~ +\,    \sum_{k\in\,\boldsymbol{{\cal I}}_m}  P\bigg\{ \widetilde{\psi}_m({\bf Z}_k; \varphi) < \varrho_0/2  \bigg\} \nonumber\\
&& ~\leq~  \exp\left\{  \frac{-2\ell^2\,(\beta/2)^2}{2L^2 \pi^{-1}_{\mbox{\tiny min}}  
\left\{ 4+B^2 (2/\varrho_0)^2  \right\}}\right\} \,+ \sum_{k\in\,\boldsymbol{{\cal I}}_m}  P\Big\{ \widetilde{\psi}_m({\bf Z}_k; \varphi) < \varrho_0/2  \Big\}. \label{ESP2}
\end{eqnarray}
If we put $\psi({\bf Z}_k; \varphi):=E\big[\Delta_k \varphi(Y_k)\big|{\bf Z}_k\big]$, then we find $ P\big\{ \widetilde{\psi}_m({\bf Z}_k; \varphi) < \varrho_0/2 \,\big| {\bf Z}_k\big\} \leq P\big\{-\widetilde{\psi}_m({\bf Z}_k; \varphi)  + \psi({\bf Z}_k; \varphi)>  \varrho_0 -\varrho_0/2  \, \big| {\bf Z}_k\big\}\leq$ $P\big\{ \big|\widetilde{\psi}_m({\bf Z}_k; \varphi)  - \psi({\bf Z}_k; \varphi)\big| > \varrho_0/2  \, \big| {\bf Z}_k\big\}\leq$ $C_{16}\exp\{-C_{17} m h^d\}$, for $n$ large enough and positive constants $C_{16}$ and $C_{17}$  not depending on $n$, where the exponential bound follows for $m$ large enough, under assumptions (B)--(E). Thus, in view of (\ref{ESP1}) and (\ref{ESP2}), one finds
\begin{eqnarray}
T_n(1) &\leq& 2\,\big|\mathcal{F}_{\varepsilon}\big|\,\left(  \exp\left\{  \frac{-2\ell^2\,(\beta/2)^2}{2L^2 \pi^{-1}_{\mbox{\tiny min}}  
\left\{ 4+B^2 (2/\varrho_0)^2  \right\}}\right\} + C_{16}\,m\, \exp\left\{-C_{17} m h^d\right\}\right). \label{Tn1t}
\end{eqnarray} 
As for the term $T_n(2)$ that appears in (\ref{Tn12}), one can use the fact that $|\widehat{m}^{\mbox{\tiny HT}}_m({\bf X}_i; \widetilde{\pi}_{\varphi}) - Y_i|^2\leq 2L^2\big\{ 4+B^2 \max^2_{k\in \boldsymbol{{\cal I}}_m}\big|1/\widetilde{\psi}_m({\bf Z}_k; \varphi)\big|   \big\} $ to write
\begin{eqnarray*}
&& T_n(2) \\
&& \leq\,  \big|\mathcal{F}_{\varepsilon}\big|\sup_{\varphi \in \mathcal{F}_{\varepsilon}} \, 
\left(P\left\{\Bigg[ \frac{2 L^2}{\ell}\Big\{ 4+
\Big(\max_{k\in \boldsymbol{{\cal I}}_m}\Big| B/\widetilde{\psi}_m({\bf Z}_k; \varphi)\Big|\Big)^2   \Big\} \,  \sum_{i\in\,\boldsymbol{{\cal I}}_\ell}
\left| \frac{1}{\pi_{\varphi}({\bf Z}_i, Y_i)}-\frac{1}{\widetilde{\pi}_{\varphi}({\bf Z}_i, Y_i)} \right| >\frac{\beta}{2}\Bigg]\right. \right.\\[2pt]
&&~~~~~~~~~~~~~~~~~~~~~~~~~~~~~~~ \left. \left. \cap \Bigg[ \bigcap_{k\in\,\boldsymbol{{\cal I}}_m}\left\{ \widetilde{\psi}_m({\bf Z}_k; \varphi) \geq\varrho_0/2 \right\}   \Bigg]
\right\} + \sum_{k\in\,\boldsymbol{{\cal I}}_m}  P\Big\{ \widetilde{\psi}_m({\bf Z}_k; \varphi) < \varrho_0/2  \Big\}  \right)\\
&& \leq\,  \big|\mathcal{F}_{\varepsilon}\big|\sup_{\varphi \in \mathcal{F}_{\varepsilon}} \, \left(P\left\{ \ell^{-1}  \sum_{i\in\,\boldsymbol{{\cal I}}_\ell}
\left| \frac{1}{\pi_{\varphi}({\bf Z}_i, Y_i)}-\frac{1}{\widetilde{\pi}_{\varphi}({\bf Z}_i, Y_i)} \right| >\frac{\beta}{4L^2[4+B^2(2/\varrho_0)^2]}\right\} 
\right.\\[2pt]
&& ~~~~~~~~~~~~~~~~~~~~~~~~~\left. + \sum_{k\in\,\boldsymbol{{\cal I}}_m}  P\Big\{ \widetilde{\psi}_m({\bf Z}_k; \varphi) < \varrho_0/2  \Big\}  \right).
\end{eqnarray*}
Employing the arguments that were used in (\ref{EQ7}), (\ref{EEQQ1}), (\ref{EEQQ2}), and (\ref{EEQQ3}), one arrives at 
\begin{equation}
T_n(2)  \,\leq\,  \big|\mathcal{F}_{\varepsilon}\big| \left( C_{46}\,\ell \exp\left\{-C_{47}\, m h^d \beta^2\right\}
+  C_{48}\,\ell \exp\left\{-C_{49}  m h^d\right\}  +  C_{50}\,\ell \exp\left\{-C_{51} m h^d \right\}
\right).   \label{Tn2t}
\end{equation} 
Now, putting together (\ref{Tn12}), (\ref{Tn1t}), and (\ref{Tn2t}), we find that for every $\beta>0$ and $n$ large enough (and thus $m$ and $\ell$),
\begin{eqnarray}
&& P\left\{ \sup_{\varphi \in \mathcal{F}_{\varepsilon}} \,\left| \widetilde{L}_{m,\ell}(\varphi) - E\left[\Big|\widehat{m}^{\mbox{\tiny HT}}_m({\bf X}; \widetilde{\pi}_{\varphi}) - Y\Big|^2\Big|\mathbb{D}_m\right] \right|  > \beta \right\} \nonumber\\
&&~~\leq\, \big|\mathcal{F}_{\varepsilon}\big| \left(  \exp\left\{-C_{52}\, \ell^2  \beta^2\right\}
+  C_{46}\,\ell \exp\left\{-C_{47}\,  m h^d\beta^2\right\}  +  C_{53}\,(\ell\vee m) \exp\left\{-C_{54}\, m h^d \right\}
\right).~~~~~~   \label{TTn12}
\end{eqnarray}
To wrap up the proof, we also need to deal with the last probability statement on the right side of (\ref{Bound99}). To that end, define the quantities
\begin{eqnarray}
Q_{n,1}(\varphi) &=& \Bigg|\frac{1}{\ell}\, \sum_{i\in\,\boldsymbol{{\cal I}}_\ell} \frac{\Delta_i \,\big|\widehat{m}^{\mbox{\tiny HT}}_{m}({\bf X}_i;\widetilde{\pi}_\varphi)-Y_i\big|^2}{ \pi_{\varphi}({\bf Z}_i, Y_i)} - \frac{1}{\ell}\, \sum_{i\in\,\boldsymbol{{\cal I}}_\ell} \frac{\Delta_i \,\big|m({\bf X}_i;\pi_\varphi)-Y_i\big|^2}{ \pi_{\varphi}({\bf Z}_i, Y_i)}\Bigg| \label{VN1}\\
Q_{n,2}(\varphi)&=&\Bigg| \frac{1}{\ell}\, \sum_{i\in\,\boldsymbol{{\cal I}}_\ell} \frac{\Delta_i \,\big|m({\bf X}_i;\pi_\varphi)-Y_i\big|^2}{ \pi_{\varphi}({\bf Z}_i, Y_i)}- E\left[\frac{\Delta\big|m({\bf X}; \pi_\varphi) - Y\big|^2}{\pi_{\varphi}({\bf Z}, Y)}\right]\Bigg| \label{VN2}\\
Q_{n,3}(\varphi) &=&\Bigg|\frac{1}{\ell} \sum_{i\in\,\boldsymbol{{\cal I}}_\ell} \Delta_i \,\big|\widehat{m}^{\mbox{\tiny HT}}_{m}({\bf X}_i;\widetilde{\pi}_\varphi)-Y_i\big|^2 \left[ \Big(\pi_{\varphi}({\bf Z}_i, Y_i)\Big)^{-1}-\Big(\widetilde{\pi}_{\varphi}({\bf Z}_i, Y_i)\Big)^{-1}\right]\Bigg|,
\end{eqnarray}
and observe that for every $\beta>0$,
\begin{eqnarray} 
&& P\left\{\sup_{\varphi \in \mathcal{F}_{\varepsilon}} \,\left| \widetilde{L}_{m,\ell}(\widetilde{\pi}_{\varphi}) - E\Big|m({\bf X}; \pi_{\varphi}) - Y\Big|^2 \right|>\beta\right\} \nonumber\\
&&  ~\leq~  P\bigg\{\sup_{\varphi \in \mathcal{F}_{\varepsilon}} \,\big| Q_{n,1}(\varphi)\big| > \frac{\beta}{3}\bigg\} 
+  P\bigg\{\sup_{\varphi \in \mathcal{F}_{\varepsilon}} \,\big| Q_{n,2}(\varphi)\big| > \frac{\beta}{3}\bigg\} 
+  P\bigg\{\sup_{\varphi \in \mathcal{F}_{\varepsilon}} \,\big| Q_{n,3}(\varphi)\big| > \frac{\beta}{3}\bigg\}.\nonumber\\[3pt]
&&~:=~P_{n,1} +P_{n,2} + P_{n,3}. ~~~~~~~\label{ABC} 
\end{eqnarray}
However, in view of (\ref{OABC}) and the fact that $|m({\bf X}_i; \pi_\varphi)|\leq L/\pi_{\mbox{\tiny $\min$}}$, one obtains
\begin{eqnarray}
&& P_{n,1} \nonumber\\
&&\leq 
P\bigg\{\sup_{\varphi \in \mathcal{F}_{\varepsilon}}\frac{1}{\ell} \sum_{i\in\boldsymbol{{\cal I}}_\ell} \Big[\Delta_i\big|\widehat{m}^{\mbox{\tiny HT}}_{m}({\bf X}_i;\widetilde{\pi}_\varphi)-m({\bf X}_i; \pi_\varphi)\big| \big|\widehat{m}^{\mbox{\tiny HT}}_{m}({\bf X}_i;\widetilde{\pi}_\varphi)+m({\bf X}_i; \pi_\varphi)\mbox{$-$\,}2Y_i\big|\Big] > \frac{\beta \pi_{\mbox{\tiny min}}}{3}\bigg\}\nonumber\\
&&\leq P\left\{\sup_{\varphi \in \mathcal{F}_{\varepsilon}}\left| \frac{1}{\ell} \sum_{i\in\boldsymbol{{\cal I}}_\ell} \left[\big| \widehat{m}^{\mbox{\tiny HT}}_{m}({\bf X}_i;\widetilde{\pi}_\varphi)-m({\bf X}_i; \pi_\varphi)\big| 
\left(3 + \frac{1}{\pi_{\mbox{\tiny min}}}+\max_{k\in \boldsymbol{{\cal I}}_m}\left|\frac{B}{\widetilde{\psi}_m({\bf Z}_k; \varphi)}\right|\right)\right]\right| > \frac{\beta \pi_{\mbox{\tiny min}}}{3L}\right\}\nonumber\\
&&\leq \big|\mathcal{F}_{\varepsilon}\big| \sup_{\varphi \in \mathcal{F}_{\varepsilon}} \sum_{i\in\boldsymbol{{\cal I}}_\ell}\Bigg(
P\Bigg\{\bigg[\Big|\widehat{m}^{\mbox{\tiny HT}}_{m}({\bf X}_i;\widetilde{\pi}_\varphi)-m({\bf X}_i; \pi_\varphi)\Big| \Big(3 + \pi_{\mbox{\tiny min}}^{-1}+B /(\varrho_0/2)\Big) > \frac{\beta \pi_{\mbox{\tiny min}}}{3L}\bigg]\nonumber\\
&&~~~~~~~~~~~~~~~~~~~~~~~~~~~~~~~~ \cap \bigg[ \bigcap_{k\in\,\boldsymbol{{\cal I}}_m}\left\{ \widetilde{\psi}_m({\bf Z}_k; \varphi) 
\geq \frac{\varrho_0}{2} \right\}   \bigg]\Bigg\} 
+ \sum_{k\in\,\boldsymbol{{\cal I}}_m}  P\Big\{ \widetilde{\psi}_m({\bf Z}_k; \varphi) < \frac{\varrho_0}{2}  \Big\}\Bigg)  \nonumber\\
&&\leq \big|\mathcal{F}_{\varepsilon}\big|\bigg[  \sup_{\varphi \in \mathcal{F}_{\varepsilon}} \sum_{i\in\boldsymbol{{\cal I}}_\ell}
P\left\{\Big|\widehat{m}^{\mbox{\tiny HT}}_{m}({\bf X}_i;\widetilde{\pi}_\varphi)-m({\bf X}_i; \pi_\varphi)\Big| > C_{\beta}\right\}
+ \ell \sum_{k\in\boldsymbol{{\cal I}}_m}  P\Big\{ \widetilde{\psi}_m({\bf Z}_k; \varphi) < \frac{\varrho_0}{2}  \Big\}\bigg]~~~~~~~~~~ \label{EQtil5}
\end{eqnarray}
where
\begin{equation*} \label{CBET}
C_{\beta} = \pi_{\mbox{\tiny min}}\beta \big/ 3L\left(3+\pi_{\mbox{\tiny min}}^{-1}+B /(\varrho_0/2)\right).
\end{equation*}
But the first probability statement in (\ref{EQtil5}) can be bounded as follows. First, observe that
\begin{eqnarray}
&& P\left\{\Big|\widehat{m}^{\mbox{\tiny HT}}_{m}({\bf X}_i;\widetilde{\pi}_\varphi)-m({\bf X}_i; \pi_\varphi)\Big| > C_{\beta}\right\}\nonumber\\
&& \leq P\left\{\Big|\widehat{m}^{\mbox{\tiny HT}}_{m}({\bf X}_i;\widetilde{\pi}_\varphi) - \widehat{m}^{\mbox{\tiny HT}}_m({\bf X}_i; \pi_\varphi)\Big| > \frac{C_{\beta}}{2}\right\}
+ P\left\{\Big|\widehat{m}^{\mbox{\tiny HT}}_{m}({\bf X}_i; \pi_\varphi)-m({\bf X}_i; \pi_\varphi)\Big| > \frac{C_{\beta}}{2}\right\}\nonumber\\
&& := \mathcal{P}_{n1}(\beta) + \mathcal{P}_{n1}(\beta). \label{EQtil6}
\end{eqnarray}
On the other hand,
\begin{eqnarray*}
\mathcal{P}_{n1}(\beta) &=& P\left\{\left|\sum_{k\in\,\boldsymbol{{\cal I}}_m}   \left[\Big(\widetilde{\pi}_{\varphi}({\bf Z}_k, Y_k)\Big)^{-1} - \Big(\pi_{\varphi}({\bf Z}_k, Y_k)\Big)^{-1} \right] \frac{\Delta_K Y_k \mathcal{K}(({\bf X}_i-{\bf X}_k)/h)}{\sum_{j\in\,\boldsymbol{{\cal I}}_m}\mathcal{K}(({\bf X}_i-{\bf X}_j)/h)  }\right|>\,\frac{C_{\beta}}{2}\right\}\\
&\leq& P\left\{\max_{k\in\,\boldsymbol{{\cal I}}_m}\left|\Big(\widetilde{\pi}_{\varphi}({\bf Z}_k, Y_k)\Big)^{-1} - \Big(\pi_{\varphi}({\bf Z}_k, Y_k)\Big)^{-1} \right|>\,\frac{C_{\beta}}{2L}\right\}\\
&\leq& \sum_{k\in\,\boldsymbol{{\cal I}}_m}  P\left\{\left|\Big(\widetilde{\pi}_{\varphi}({\bf Z}_k, Y_k)\Big)^{-1} - \Big(\pi_{\varphi}({\bf Z}_k, Y_k)\Big)^{-1} \right|>\,\frac{C_{\beta}}{2L}\right\}.
\end{eqnarray*}
Therefore, using arguments similar to those leading to (\ref{EQ7}), (\ref{EEQQ1}), (\ref{EEQQ2}), and (\ref{EEQQ3}), we find, for every $\beta>0$ and $n$ large enough,
\begin{equation*}
\mathcal{P}_{n1}(\beta)\, \leq\, C_{39}\,  m e^{-C_{40} m h^d \beta^2} +C_{41}\, m e^{-C_{42}m h^d},
\end{equation*}
where $C_{39}$--\,$C_{42}$ are positive constants not depending on $m$, $\ell$, or $\beta$. Furthermore, tedious but standard arguments can be used to show that for $n$ large enough, there are positive constants $C_{43}$ and $C_{44}$,  not depending on $m$, $\ell$, or $\beta$, such that
\begin{equation*}
\mathcal{P}_{n2}(\beta)\, \leq\, C_{43}\, e^{-C_{44} m h^d \beta^2}.
\end{equation*}
As for the last probability statement on the right side of (\ref{EQtil5}), our earlier arguments (see the  paragraph after equation (\ref{ESP2})) yield
$P\big\{ \widetilde{\psi}_m({\bf Z}_k; \varphi) < \varrho_0/2  \big\} \leq C_{16}\,\exp\{-C_{17}\, m h^d\}$, for $n$ large enough, where $C_{16}$ and $C_{17}$ are positive constants not depending on $n$. Therefore, in view of (\ref{EQtil5}) we arrive at 
\begin{equation} \label{EQtil7}
P_{n,1} \,\leq\, \ell\, \big|\mathcal{F}_{\varepsilon}\big|  \left( C_{39}\,m e^{-C_{40} m h^d \beta^2 } + C_{43}\, e^{-C_{44} m h^d \beta^2 } +   C_{55}\,m e^{-C_{56} m h^d  } \right),
\end{equation}
for $n$ large enough, where $P_{n,1}$ is as in (\ref{ABC}). To deal with $P_{n,2}$, we first note that the terms $\Delta_i \,\big|m({\bf X}_i;\pi_\varphi)-Y_i\big|^2/ \pi_{\varphi}({\bf Z}_i, Y_i)$, $i\in\,\boldsymbol{{\cal I}}_\ell$, are iid bounded random variables taking values in the interval $\big[0,\,L^2(1+1/\pi_{\mbox{\tiny min}})^2/\pi_{\mbox{\tiny min}}\big]$. Therefore an application of Hoeffding's inequality (in conjunction with the union bound) immediately yields
\begin{eqnarray} 
P_{n,2}  &\leq&  \big|\mathcal{F}_{\varepsilon}\big| \sup_{\varphi \in \mathcal{F}_{\varepsilon}}
P\left\{\bigg|\frac{1}{\ell}  \sum_{i\in\boldsymbol{{\cal I}}_\ell}  \frac{\Delta_i \,\big|m({\bf X}_i;\pi_\varphi)-Y_i\big|^2}{ \pi_{\varphi}({\bf Z}_i, Y_i)}- E\left[\frac{\Delta\big|m({\bf X}; \pi_\varphi) - Y\big|^2}{\pi_{\varphi}({\bf Z}, Y)}\right]\Bigg| > \frac{\beta}{3}\right\} \nonumber\\
&\leq& 2\,\big|\mathcal{F}_{\varepsilon}\big|  \exp\Big\{-2\pi_{\mbox{\tiny min}}^2 \ell\, (\beta/3)^2\big/[L^4 (1+1/\pi_{\mbox{\tiny min}})^4]\Big\}. \label{ADC}
\end{eqnarray}
Finally, to deal with the term $P_{n,3}$ in (\ref{ABC}), we observe that  in view of (\ref{OABC}), and with $\varrho_0$ as in Assumption (D), one has
\begin{eqnarray*}
P_{n,3}  &\leq &  \big|\mathcal{F}_{\varepsilon}\big| \sup_{\varphi \in \mathcal{F}_{\varepsilon}}
\left(P\left\{\Bigg[ \bigg(2 +\max_{k\in \boldsymbol{{\cal I}}_m}\bigg|\frac{B}{\widetilde{\psi}_m({\bf Z}_k; \varphi)}\bigg| \bigg)^2 L^2 \cdot \frac{1}{\ell} \sum_{i\in \boldsymbol{{\cal I}}_\ell}\bigg|\frac{1}{\widetilde{\pi}_{\varphi}({\bf Z}_i, Y_i)} - \frac{1}{\pi_{\varphi}({\bf Z}_i, Y_i)} \bigg| > \frac{\beta}{3}\Bigg]\right .\right. \nonumber\\
&&~~~~~~~~~~~~~~~~~~~~\left.\left. \cap \Bigg[ \bigcap_{k\in\,\boldsymbol{{\cal I}}_m}\left\{ \widetilde{\psi}_m({\bf Z}_k; \varphi) \geq\varrho_0/2 \right\}   \Bigg] \right\} + \sum_{k\in\,\boldsymbol{{\cal I}}_m}  P\Big\{ \widetilde{\psi}_m({\bf Z}_k; \varphi) < \varrho_0/2  \Big\}\right)  \nonumber\\
&\leq &  \big|\mathcal{F}_{\varepsilon}\big| \sup_{\varphi \in \mathcal{F}_{\varepsilon}}
\left(P\left\{ \frac{1}{\ell} \sum_{i\in \boldsymbol{{\cal I}}_\ell}   \bigg|\frac{1}{\widetilde{\pi}_{\varphi}({\bf Z}_i, Y_i)} - \frac{1}{\pi_{\varphi}({\bf Z}_i, Y_i)} \bigg| > d_{\beta}\right\} + \sum_{k\in\boldsymbol{{\cal I}}_m}  P\Big\{ \widetilde{\psi}_m({\bf Z}_k; \varphi) < \varrho_0/2  \Big\}\right)
\end{eqnarray*}
where $d_{\beta} = \big[ 3L^2 (2 +2 B/\varrho_0)^2  \big]^{-1}\beta$. Now, employing the arguments used to bound the term $S_n(2)$ in (\ref{EQ5}), (see (\ref{EQ7}), (\ref{EEQQ1}), (\ref{EEQQ2}), (\ref{EEQQ3})), it is straightforward to show that for $n$ large enough
\begin{equation}
P_{n,3} \,\leq\,  \big|\mathcal{F}_{\varepsilon}\big| \Big( C_{58} \ell \,e^{-C_{59} m h^d \beta^2} + C_{60} m\, e^{-C_{61} m h^d} \Big), \label{EQtil8}
\end{equation}
for positive constants $C_{58}$--$\,C_{61}$ not depending on $\ell$, $m$, or $\beta$. Putting together (\ref{ABC}), (\ref{EQtil7}), (\ref{ADC}), and (\ref{EQtil8}), one finds that for each  $\beta>0$ and $n$ large enough, 
\begin{eqnarray}
P\left\{
\sup_{\varphi \in \mathcal{F}_{\varepsilon}} \,\left| \widetilde{L}_{m,\ell}(\widetilde{\pi}_{\varphi}) - 
E\Big|m({\bf X}; \pi_{\varphi}) - Y\Big|^2 \right|>\beta\right\} &\leq& \big|\mathcal{F}_{\varepsilon}\big| \Big( C_{39} \ell m \,e^{-C_{40} m h^d \beta^2} + C_{55} \ell m\, e^{-C_{56} m h^d} \nonumber\\
&&~~~~~~~~~~+\,  2 \,e^{-C_{64} \ell \beta^2}\Big). \label{NOP1}
\end{eqnarray}
Now, for any decreasing sequence $0<\varepsilon_n \downarrow 0$,   let $\varphi_{\varepsilon_n}$ be as in (\ref{FIEP}). Then, employing arguments similar to those used to arrive at    (\ref{EQ4B}) and (\ref{EEQ2}), give
 \begin{eqnarray}
\int\Big|\widehat{m}^{\mbox{\tiny HT}}({\bf x}; \widetilde{\varphi}_n)- m({\bf x}) \Big|^2  \mu(d{\bf x}) &=&
\int\Big|\widehat{m}^{\mbox{\tiny HT}}({\bf x}; \widetilde{\varphi}_n) - m({\bf x}; \varphi_{\varepsilon_n}) +  m({\bf x}; \varphi_{\varepsilon_n})- m({\bf x}) \Big|^2  \mu(d{\bf x}) \nonumber\\
&\leq& 2\int\Big|\widehat{m}^{\mbox{\tiny HT}}({\bf x}; \widetilde{\varphi}_n) - m({\bf x}; \varphi_{\varepsilon_n}) \Big|^2  \mu(d{\bf x})\,+ \,4LC \,\varepsilon_n\,, \label{STARSTAR}
\end{eqnarray}
where $C>0$ is the constant in Lemma \ref{LEM-2}.  Therefore, in view of (\ref{STARSTAR}) and (\ref{Bound99}), for every constant $t>0$ we have
 \begin{eqnarray*}
&& \frac{1}{2}\,
P\left\{\int\Big| \widehat{m}^{\mbox{\tiny HT}}({\bf x}; \widetilde{\pi}_{\widetilde{\varphi}_{n}}) - m({\bf x}) \Big|^2  \mu(d{\bf x}) \,>\,  t \right\}\\
&& ~~~~~\leq~ \frac{1}{2}\, P\left\{\int\Big| \widehat{m}^{\mbox{\tiny HT}}({\bf x}; \widetilde{\pi}_{\widetilde{\varphi}_{n}}) - m({\bf x}; \varphi_{\varepsilon_n}) \Big|^2  \mu(d{\bf x})  >\, t/2-2LC\varepsilon_n \right\}\\
&&~~~~~\leq~ P\left\{\int\Big| \widehat{m}^{\mbox{\tiny HT}}({\bf x}; \widetilde{\pi}_{\widetilde{\varphi}_{n}}) - m({\bf x}; \varphi_{\varepsilon_n}) \Big|^2  \mu(d{\bf x})  >\, t/2-2LC\varepsilon_n \right\}\\
&& ~~~~~~~~~~~~~~ - P\left\{\int\Big| \widehat{m}^{\mbox{\tiny HT}}({\bf x}; \widetilde{\pi}_{\widetilde{\varphi}_{n}}) - m({\bf x}; \varphi_{\varepsilon_n}) \Big|^2  \mu(d{\bf x})  >\, (t/2-2LC\varepsilon_n)^2/(c_4 \,\varepsilon_n) \right\}~~~~~~~~\\
&&~~~~~~~~~~~~~~~~(\mbox{for $n$ large enough, where $c_4>0$ is as in the first line of  (\ref{Bound99})})\\
&&~~~~~\leq~ P\left\{ \sup_{\varphi \in \mathcal{F}_{\varepsilon_n}} \,\left|E\left[\Big|\widehat{m}^{\mbox{\tiny HT}}_m({\bf X}; \widetilde{\pi}_{\varphi}) - Y\Big|^2\Big|\mathbb{D}_m\right] - \widetilde{L}_{m,\ell}(\widetilde{\pi}_{\varphi}) \right|>\frac{t/2-2LC\varepsilon_n}{3}\right\}~~~~~\nonumber \\[3pt]
&&~~~~~~~~~~~~~~+   
P\left\{\sup_{\varphi \in \mathcal{F}_{\varepsilon_n}} \,\left| \widetilde{L}_{m,\ell}(\widetilde{\pi}_{\varphi}) - 
E\Big|m({\bf X}; \pi_{\varphi}) - Y\Big|^2 \right|>\frac{t/2-2LC\varepsilon_n}{3}\right\}.
\end{eqnarray*}
Finally, choosing $n$ large enough so that $(t/2-2LC\varepsilon_n)/3 > t/12$, and using the bounds in (\ref{NOP1}) and (\ref{TTn12}), we find
\begin{eqnarray*}
P\left\{\int\Big| \widehat{m}^{\mbox{\tiny HT}}({\bf x}; \widetilde{\pi}_{\widetilde{\varphi}_{n}}) - m({\bf x}) \Big|^2  \mu(d{\bf x}) >\,  t \right\}
&\leq& \big|\mathcal{F}_{\varepsilon_n}\big| \Big( C_{65}\,e^{-C_{66} \ell t^2} + C_{67} e^{-C_{68} \ell^2 t^2}+ C_{69} \ell \,e^{-C_{70}m h^d t^2} \\
&& ~ + C_{71}\,\ell m \,e^{-C_{72} m h^d (t^2\vee 1)}  + C_{73} (\ell \vee m) e^{-C_{74} m h^d}\Big),
\end{eqnarray*}
for $n$ large enough where $C_{65}$--$\,C_{74}$ are positive constants not depending on $m$ , $\ell$, or $t$. This completes the proof of Part (i) of the theorem.

\vspace{5mm}\noindent
{\it Part (ii).}

\vspace{1mm}\noindent
The proof of Part (ii) of the theorem is virtually the same and, in fact, easier and therefore will not be given.

\hfill $\Box$

\vspace{4mm}\noindent
PROOF OF COROLLARY \ref{COR-AB}

\vspace{1.5mm}\noindent
The corollary follows from the Borel-Cantelli lemma in conjunction with (\ref{cond2}), the bound in Theorem \ref{THM-BBC}, and Remark \ref{REM-thm4}.

\hfill $\Box$

\vspace{4mm}\noindent
PROOF OF THEOREM \ref{THM-BBCC}

\vspace{1.5mm}\noindent
The proof of this theorem is similar to that of Theorem \ref{THM-BB} and therefore will not be given.

\hfill $\Box$

\vspace{4mm}\noindent
PROOF OF THEOREM \ref{THM-class-1}

\vspace{1.5mm}\noindent
{\it Part (i).}

\vspace{1mm}\noindent
By (\ref{Bayes-bound}), we have
\begin{equation} \label{p2bound}
P\left\{\widehat{g}_n({\bf X};\widehat{\varphi}_{n}) \neq Y\Big| \mathbb{D}_n \right\} - P\{g_{\mbox{\tiny B}}({\bf X})\neq Y\} ~\leq~ 2E\left[\Big|\widehat{m}({\bf X};\widehat{\varphi}_{n})- m(X)\Big|\,\bigg|\mathbb{D}_n \right].
\end{equation}
Now, Part (i) of the theorem follows from (\ref{p2bound}) and  Corollary \ref{COR-A} in conjunction with the Cauchy-Schwarz inequality.

\vspace{4mm}\noindent
{\it Part (ii).} 

\vspace{1mm}\noindent
Taking the expectation of both sides of (\ref{p2bound}), the result follows from Corollary \ref{COR-A2} together with the Cauchy-Schwarz inequality.

\vspace{4mm}\noindent
{\it Part (iii).} 

\vspace{1mm}\noindent
By a result of Audibert and Tsybakov (2007; Lemma 5.2), under the margin assumption (G), we have 
\begin{equation} \label{p3bound}
P\left\{\widehat{g}_n({\bf X};\widehat{\varphi}_{n}) \neq Y \right\} - P\{g_{\mbox{\tiny B}}({\bf X})\neq Y\} ~\leq~ \left(E\Big|\widehat{m}({\bf X};\widehat{\varphi}_{n})- m(X)\Big|^2\right)^{\frac{1+\alpha}{2+\alpha}},
\end{equation}
where $\alpha$ is as in (\ref{MARG}). The result now follows from Corollary \ref{COR-A2}.

\hfill     $\Box$

\vspace{4mm}\noindent
PROOF OF THEOREM \ref{THM-class-2}

\vspace{1.5mm}\noindent
The proof uses Corollaries \ref{COR-AB} and \ref{COR-AB2} and is virtually the same as that of Theorem \ref{THM-class-2}, and thus will not be given.

\hfill     
$\Box$

\vspace{6mm}\noindent
{\bf \large Appendix.} 

\vspace{1mm}\noindent
PROOF OF (\ref{eps-cover})\\
To show that (\ref{eps-cover}) is an $\varepsilon$-cover of the class $\mathcal{F}$ in (\ref{Exam-1}), let  
\[
\Omega_{\varepsilon} = \left\{
2\,i \varepsilon/(L\exp(ML))\,\bigg|
- \floor*{ML\exp(ML)/\varepsilon}\,\leq\, i\,\leq \floor*{ML\exp(ML)/\varepsilon}    \right\}\cup\,\{-M\}\,\cup\,\{M\}.
\]
Also, let
$\gamma\in[-M , M]$ be given and put $\varphi(y) = e^{\gamma y} \in \mathcal{F}$. If $\widetilde{\gamma}\in \Omega_{\varepsilon}$ is the closest value to $\gamma$, then 
\begin{eqnarray*}
	\sup_{ |y| \leq L} \Big|e^{\gamma y}-e^{\widetilde{\gamma}y} \Big|
	&=& \sup_{|y| \leq L}  \Big| y\, \exp\{\gamma^{\dagger} y\} \Big|\cdot\big| \widetilde{\gamma}-\gamma\big|,~~~
	\mbox{where}~~\gamma^{\dagger}\in(\widetilde{\gamma}\wedge \gamma\,,\,   \widetilde{\gamma}\vee \gamma)\\[2pt]
	&\leq& L\exp\{ML\}\cdot \big|\widetilde{\gamma}-\gamma\big|\\
	&\leq& L\exp\{ML\}\cdot \frac{\varepsilon}{L\exp\{ML\}} ~=~\varepsilon\,,\\
\end{eqnarray*}
where the last line follows from the fact that the distance between $\gamma$ and its nearest value in $\Omega_{\varepsilon}$ is bounded by $\varepsilon/(L \exp\{ML\})$.
Therefore, the class $\mathcal{F}$ is totally bounded. Moreover, a count of the number of terms in $\Omega_{\varepsilon}$ shows that the $\varepsilon$-covering number of $\mathcal{F}$ is  bounded by the quantity $2\floor*{ML\exp\{ML\}\varepsilon^{-1}}+3$.

\hfill
$\Box$

\vspace{7mm}\noindent
\noindent
{\bf Acknowledgements}\\
This work was supported by the National Science Foundation Grant DMS-1916161 of Majid Mojirsheibani.

\vspace{10mm}\noindent
{\bf Conflict of interest}\\
On behalf of all authors, the corresponding author states that there is no conflict of interest.

\vspace{12mm}\noindent
{\bf References}

\vspace{2.5mm}\noindent
Azizyan, M., Singh, A., Wasserman, L., et al. (2013) Density-sensitive semisupervised inference. \textit{Ann. Statist.} {\bf 41} 751--771.\\

\noindent
Audibert, J. Y. and  Tsybakov, A. B. (2007). Fast learning rates for plug-in classifiers under the margin condition.\textit{Ann. Statist.} \textbf{35} 608--633.\\

\noindent
Chen, X., Diao, G., and Qin, J. (2020). Pseudo likelihood-based estimation and testing of missingness mechanism function in nonignorable missing data problems. \textit{Scand. J. Stat.} \textbf{47} 1377--1400.\\

\noindent
Devroye, L., Gy\"{o}rfi, L., and Lugosi, G. (1996) A probabilistic theory of pattern recognition. Springer-Verlag, New York.\\

\noindent
Devroye, L. and Krzy\`{z}ak, A. (1989). An equivalence theorem for $L_1$ convergence of kernel regression estimate. {\it Journal of Statistical Planning and Inference,} 23, 71-82.\\

\noindent
D\"{o}ring, M., Gy\"{o}rfi, L.,  and Walk, H. Exact rate of convergence of kernel-based classification rule.  Challenges in computational statistics and data mining, 71–91,
Stud. Comput. Intell., 605, Springer, Cham, 2016.\\

\noindent
Fang, F., Zhao, J., and Shao, J. (2018). Imputation-based adjusted score equations in generalized linear models with nonignorable missing covariate values. \textit{Statistica Sinica}. \textbf{28} 1677--1701.\\

\noindent
Horvitz D. G. and  Thompson D. J. (1952). A generalization of sampling without replacement from a finite universe. {\it J. Am. Statist. Assoc}. \textbf{47} 663--685\\

\noindent
Kim, J.K. and  Yu, C.L. (2011). A semiparametric estimation of mean functionals with nonignorable missing data. {\it J. Am. Statist. Assoc}. \textbf{106} 157--65.\\

\noindent
Kohler, M. and Krzy\.{z}ak, A. (2007). On the rate of convergence of local averaging plug-in classification rules under a margin condition. \textit{IEEE Trans. Inform. Theory} \textbf{53} 1735--1742.\\

\noindent
Liu, Z. and Yau, C.-Y.  (2021). Fitting time series models for longitudinal surveys with nonignorable missing data. \textit{J. Statist. Plann. Inference.} \textbf{214} 1--12.\\

\noindent
Maity, A., Pradhan, V., and Das, U.  (2019). Bias reduction in logistic regression with missing responses when the missing data mechanism is nonignorable. \textit{Amer. Statist.} \textbf{73} 340--349.\\

\noindent
Mammen, E. and Tsybakov, A.B. (1999)\,Smooth discriminant analysis. \textit{Ann.\,Statist.}\,\textbf{27} 1808-1829.\\

\noindent
Massart, P. and E. N\'{e}d\'{e}lec, E. (2006). Risk bounds for statistical learning. \textit{Ann. Statist.} \textbf{34} 2326--2366.\\

\noindent
Mitrinovic, D. S. Analytic Inequalities. New York. Springer-Verlag, 1970.\\

\noindent
Mojirsheibani, M. (2021). On classification with nonignorable missing data. \textit{J. Multivariate Anal.} \textbf{184} 104755.  \\

\noindent
Morikawa, K., Kim, J. K., and  Kano, Y. (2017).  Semiparametric maximum likelihood estimation with data missing not at random. \textit{Can. J. Statist.} \textbf{45} 393--409.\\

\noindent
Morikawa, K. and Kim, J. K. (2018). A note on the equivalence of two semiparametric estimation methods for nonignorable nonresponse. \textit{Stat. \& Probab. Lett.} \textbf{140} 1--6.\\

\noindent
Nadaraya, E. A. (1964). On estimating regression. \textit{\it Theory Probab. Appl.} \textbf{9} 141--142.\\

\noindent
O'Brien, J., Gunawardena, H., Paulo, J., Chen, X., Ibrahim, J., Gygi, S., and Qaqish, B. (2018). The effects of nonignorable missing data on label-free mass spectrometry proteomics experiments. \textit{Ann. Appl. Statist.} \textbf{12} 2075--2095.\\

\noindent
Sadinle, M. and Reiter, J.  (2019). Sequentially additive nonignorable missing data modelling using auxiliary marginal information.  \textit{Biometrika}. \textbf{106} 889--911.\\

\noindent
Shao, J. and Wang, L. (2016) Semiparametric inverse propensity weighting for nonignorable missing data.  \textit{Biometrika.} \textbf{103} 175--187.\\

\noindent
Tsybakov, A.B. and  van de Geer, S. (2005). Square root penalty: adaptation to the margin in classification and in edge estimation. \textit{Ann. Statist.} \textbf{33} 1203--1224.\\

\noindent
Uehara, M. and Kim, J.K. (2018). Semiparametric response model with nonignorable nonresponse. Preprint on arXiv:1810.12519. ~
https://arxiv.org/abs/1810.12519v1 \\

\noindent
van der Vaart, A., Wellner, J. (1996) Weak Convergence and Empirical Processes with Applications to Statistics. Springer, New York.\\

\noindent
Watson, G.S. (1964). Smooth regression analysis. {\it Sankhya, Ser. A.} \textbf{26} 359--372. \\

\noindent
Wang, L., Shao, J., and Fang, F. (2021). Propensity model selection with nonignorable nonresponse and instrument variable.  \textit{Statistica Sinica} \textbf{31} 647--671.\\

\noindent
Wang, S., Shao, J., and Kim, J.K. (2014). Identifiability and estimation in problems with nonignorable nonresponse. Statistica Sinica 24, 1097 - 1116.\\

\noindent
Wang, J. and Shen, X. (2007) Large margin semi-supervised  learning. {\it J. Mach. Learn. Res.}, {\bf 8} 1867--1891.\\

\noindent
Yuan, C., Hedeker, D., Mermelstein, R., Xie, H.  (2020). A tractable method to account for high-dimensional nonignorable missing data in intensive longitudinal data. \textit{Stat. Med.} \textbf{39} 2589--2605.\\

\noindent
Zhao, J., Shao, J. (2015). Semiparametric pseudo-likelihoods in generalized linear models with nonignorable missing data. {\it J. Am. Statist. Assoc} 110, 1577-1590.\\

\noindent
Zhao, P., Wang, L.,  and Shao, J.  (2019). Empirical likelihood and Wilks phenomenon for data with nonignorable missing values. \textit{Scand. J. Stat.} \textbf{46}  1003--1024.

\end{document}